\documentclass{amsproc}

\usepackage{amssymb}
\usepackage{amsthm}
\usepackage{amsmath}
\usepackage{graphicx}
\usepackage{color}
\usepackage{enumerate}
 \numberwithin{equation}{section}
\usepackage{booktabs}
\allowdisplaybreaks

\usepackage[pdftex]{hyperref}
 \usepackage{mathtools}
 \mathtoolsset{showonlyrefs}
 \usepackage[nocompress]{cite}
  \usepackage{floatrow}%
  \usepackage{caption}
\usepackage{subcaption}
\usepackage{url}
\theoremstyle{plain}
\newtheorem{thm}{Theorem}[section]
\newtheorem{cor}[thm]{Corollary}

\newtheorem{lem}[thm]{Lemma}
\newtheorem{prop}[thm]{Proposition}

\theoremstyle{definition}
\newtheorem{defn}[thm]{Definition}

\theoremstyle{remark}

\newtheorem{rem}[thm]{Remark}

\newcommand{\N}{\mathbb{N}}
\newcommand{\R}{\mathbb{R}}

\newcommand{\bp}{\begin{proof}[\ensuremath{\mathbf{Proof}}]}
\newcommand{\bs}{\begin{proof}[\ensuremath{\mathbf{Solution}}]}
\newcommand{\ep}{\end{proof}}
\newcommand{\be}{\begin{equation}}
\newcommand{\ee}{\end{equation}}
\newcommand{\Stwo}{{\mathbb{S}^2}}
\newcommand{\Pee}{{\mathcal P}}
\newcommand{\id}{\text{id}_{u^\perp}}
\newcommand{\Cee}{{\mathcal C}}

\newcommand{\Jay}{{\mathcal J}}
\newcommand{\Vee}{{\mathcal V}}
\newcommand{\bg}{{\bf g}}
\newcommand{\bh}{{\boldsymbol h}}
\newcommand{\bmu}{{\boldsymbol \mu}}
\newcommand{\bxi}{{\boldsymbol \xi}}

\newcommand{\bzeta}{{\boldsymbol \zeta}}

\begin{document}

\title{A doubly monotone flow for constant width bodies in $\mathbb{R}^3$}


\author{Ryan Hynd}
\address{209 South 33rd Street Philadelphia, PA 19104-6395}
\curraddr{}
\email{rhynd@math.upenn.edu}
\thanks{}


\date{}

\begin{abstract}
We introduce a flow in the space of constant width bodies in three-dimensional Euclidean space that simultaneously increases the volume and decreases the circumradius of the shape as time increases. Starting from any initial constant width figure, we show that the flow exists for all positive times and converges to a closed ball as time tends to plus infinity.  We also anticipate that this flow is interesting to study for negative times and that it would provide a mechanism to decrease the volume and increase the circumradius of a constant width body.
\end{abstract}

\maketitle

\tableofcontents

\section{Introduction}
A constant width body is a compact, convex subset of Euclidean space in which parallel supporting planes are separated by the same distance in every direction.  In this note, we will focus on bodies of width one and simply refer to them as having constant width.  The simplest example of a constant width body is a (closed) ball of radius $1/2$. It is also known that balls of radius $1/2$ encloses the most volume of any constant width shape.

\par It is natural to inquire about volume-minimizing constant width bodies.  In the plane, it was proved independently by Lebesgue \cite{Lebesgue,Lebesgue2} and Blaschke \cite{MR1511839,Blaschke2} over a century ago that Reuleaux triangles encloses the least amount of area. A Reuleaux triangle is the intersection of  three closed disks of radius 1 which are centered at the vertices of an equilateral triangle of side length one.  There have been many subsequent proofs of the Lebesgue-Blasckhe theorem including \cite{MR1568234, MR205152,MR1391153,MR1881292,MR51543,MR2559951}. We also note that Harrell \cite{MR1881292} showed that Reuleaux triangles are uniquely area-minimizing among constant width shapes.

\begin{figure}[h]
\centering
 \includegraphics[width=.4\textwidth]{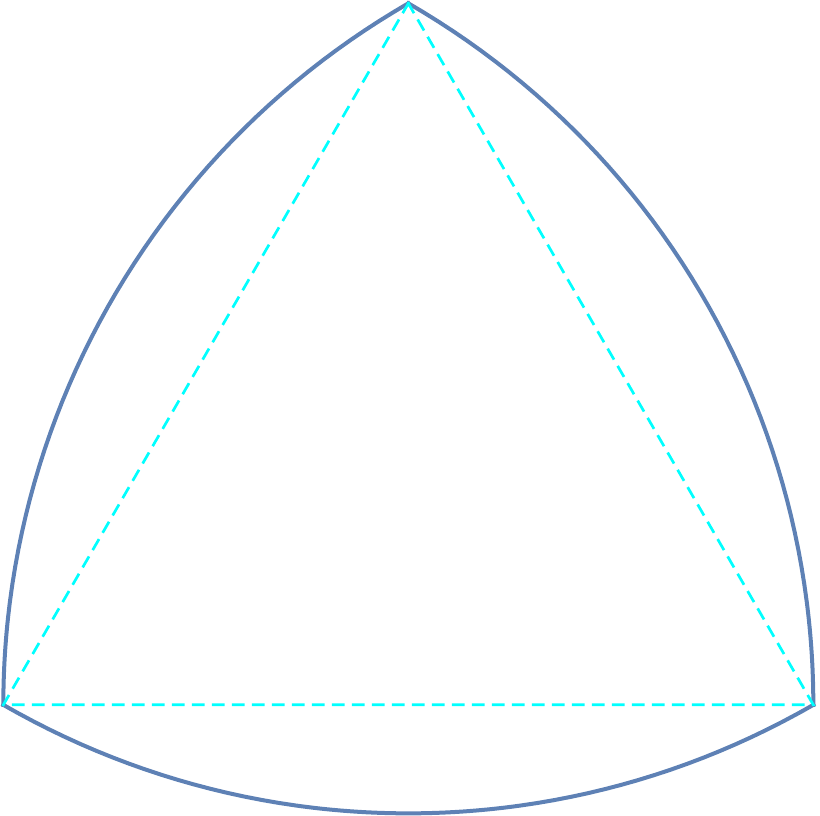}
 \caption{The boundary curve of a Reuleaux triangle with inscribed equilateral triangle.}\label{Reuleauxtriangle}
\end{figure}

\par For constant width bodies in dimensions larger than two, it is known that volume-minimizing constant width bodies exist. 
This can be seen as a consequence of the Blaschke selection theorem. While there have been some notable work on this topic such as \cite{MR2763770,MR2342202,MR1371579,MR2883369,MR205152,anciaux2009blaschkelebesgue}, surprisingly little is known about these shapes.  However, there are conjectured volume-minimizing constant width bodies in three-dimensional Euclidean space \cite{MR2844102,MR3930585,MR1316393,MR920366,MR0123962}. These are the bodies that Meissner (and Schilling) constructed  \cite{Meissner} which are based on a regular tetrahedron, somewhat analogous to how the Reuleaux triangle is based on an equilateral triangle.

\par First Meissner considered a Reuleaux tetrahedron, which is the intersection four balls of radius one each centered at the vertices of a regular tetrahedron.  This figure has four faces, four vertices and six edges just like the regular tetrahedron. After realizing the Reuleaux tetrahedron does not have constant width, Meissner was able to round three of the six edges in two ways to obtain two distinct constant width bodies; we recommend diagrams 106 and 107 of \cite{MR0123962} for a detailed description of these procedures.  Both of Meissner's tetrahedra have the same volume and surface area, and it has been long thought that these are volume-minimizing shapes. There have also been at least two numerical studies \cite{antunes2019parametric,Muller} which support this conjecture.

\par An interesting feature of a constant width body is that its inball and circumball are concentric. Moreover, the radii of these two balls sum to one.  It turns out that any body of constant width which includes 
the regular simplex of diameter one necessarily has the largest possible circumradius (and therefore the smallest possible inradius). For example, the Reuleaux triangle has this property and so do Meissner's tetrahedra. It then seems reasonable to investigate the following question.
\\\\
\noindent {\bf Question}. {\it Is there is a connection between least volume bodies of constant width and those having largest circumradius?}
\\
\par In an attempt to develop an approach to this question, we will propose a flow in the space of constant width bodies in $\R^3$ which has two distinctive features: when time moves forward, the volume increases and the circumradius decreases along the flow.  We expect that as time tends to infinity, the flow would deform any starting shape into a ball of radius 1/2. We will investigate the existence of this flow and its behavior for large times in detail below. Nevertheless, we present this flow as a possible device which can be used to answer our motivating question.  Namely, upon reversing time, we wonder if the limiting shapes exist and lead us to some insight on understanding a possible relationship between least volume bodies of constant width and those having largest circumradius.

\begin{figure}
     \centering
     \begin{subfigure}[b]{0.5\textwidth}
         \centering
         \includegraphics[width=\textwidth]{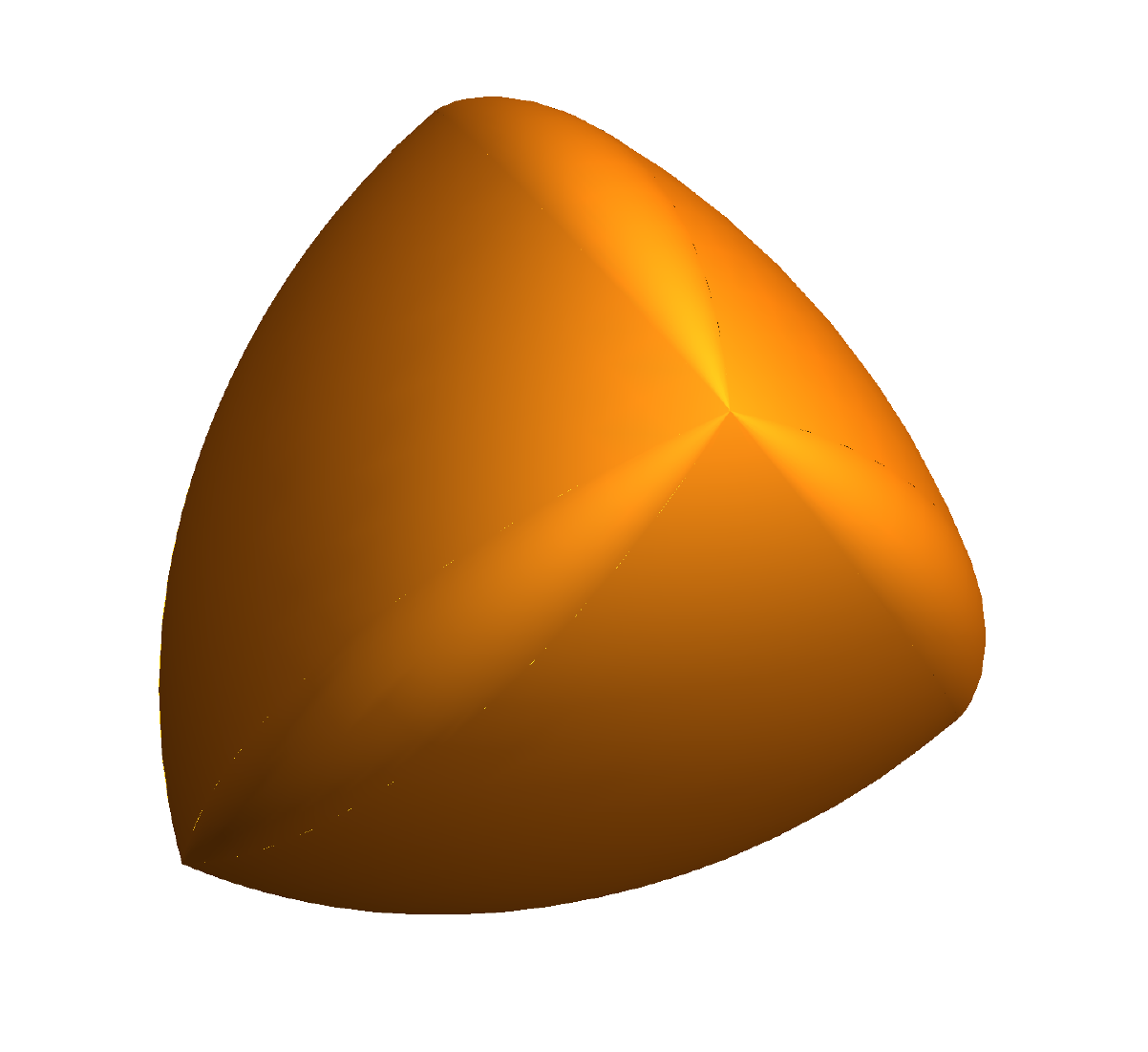}
     \end{subfigure}
     \begin{subfigure}[b]{0.49\textwidth}
         \centering
         \includegraphics[width=\textwidth]{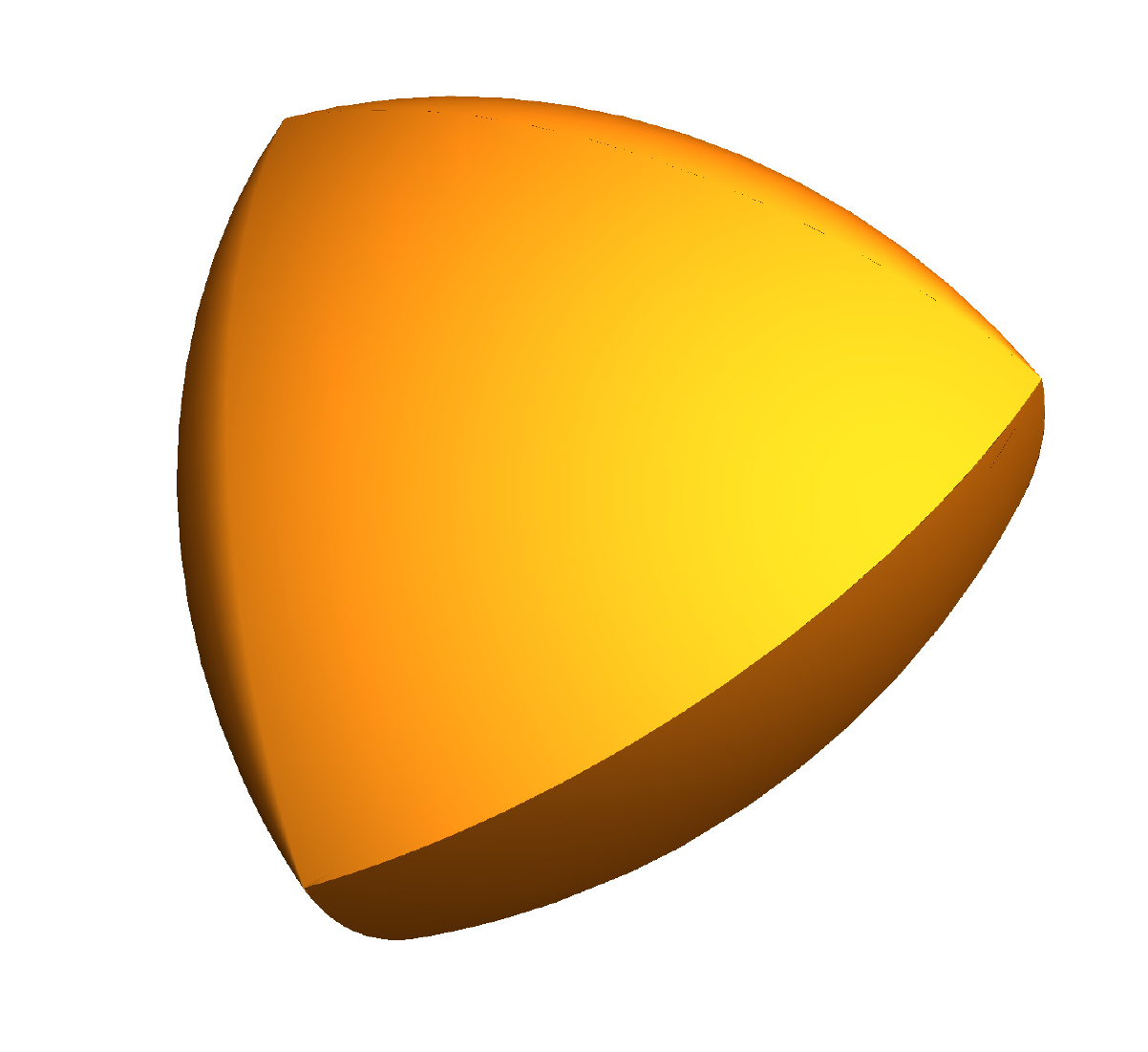}
     \end{subfigure}
          \caption{A Meissner tetrahedron in which three rounded edges meet in a vertex.}
        \label{FirstMeissnerFigure}
\end{figure}

\subsection{Quantities of interest}
In what follows, we will employ the support function 
$$
H(u)=\max_{x\in K}x\cdot u\quad (u\in \R^3) 
$$
of a constant width body $K\subset \R^3$. As $H$ is positively homogeneous, it is determined by its restriction $h=H|_{\Stwo}$. The constant width property of $K$ is equivalent to $h$ satisfying
\be
h(u)+h(-u)=1
\ee
for all $u\in\Stwo$.  It will also be convenient to consider
\be
g=h-1/2,
\ee
which is an odd function on $\Stwo$.

\par Later in this note, we will show that the circumradius of $K$ can be expressed in terms of $g$ as 
$$
R(K)=\frac{1}{2}+\min_{a\in \R^3}\max_{|u|=1}|g(u)+a\cdot u|.
$$
Likewise, we will explain that volume enclosed by $K$ can be written in terms of $g$ as
$$
V(K)=\frac{\pi}{6}-\frac{1}{2}\int_{\Stwo}\left(\frac{1}{2}|\nabla g|^2-g^2\right)d\sigma.
$$
Here $\sigma$ is two-dimensional Hausdorff measure on $\R^3$ normalized so that $\sigma(\Stwo)=4\pi$.

\par In addition, we will recall that the Hausdorff distance $d_{\mathcal H}(K_1,K_2)$ between two convex bodies $K_1,K_2\subset \R^3$ is given by the supremum norm  over $\Stwo$ of the difference of their respective support functions.  It follows 
that
$$
d_{\mathcal H}(K_1,K_2)=\max_{|u|=1}|g_1(u)-g_2(u)|.
$$
\\
\begin{figure}[h]
     \centering
     \begin{subfigure}[b]{0.5\textwidth}
         \centering
         \includegraphics[width=\textwidth]{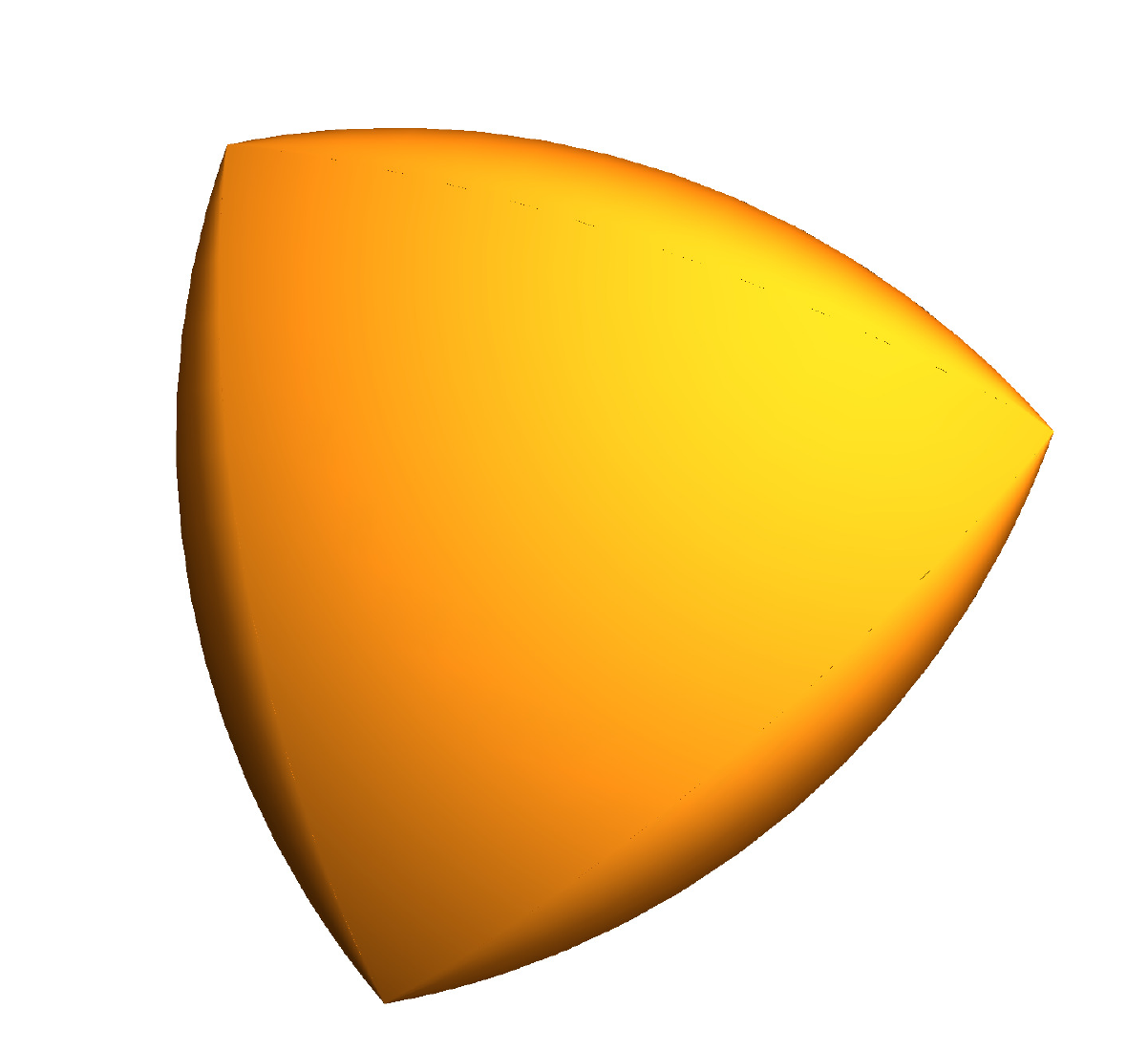}
     \end{subfigure}
     \begin{subfigure}[b]{0.49\textwidth}
         \centering
         \includegraphics[width=\textwidth]{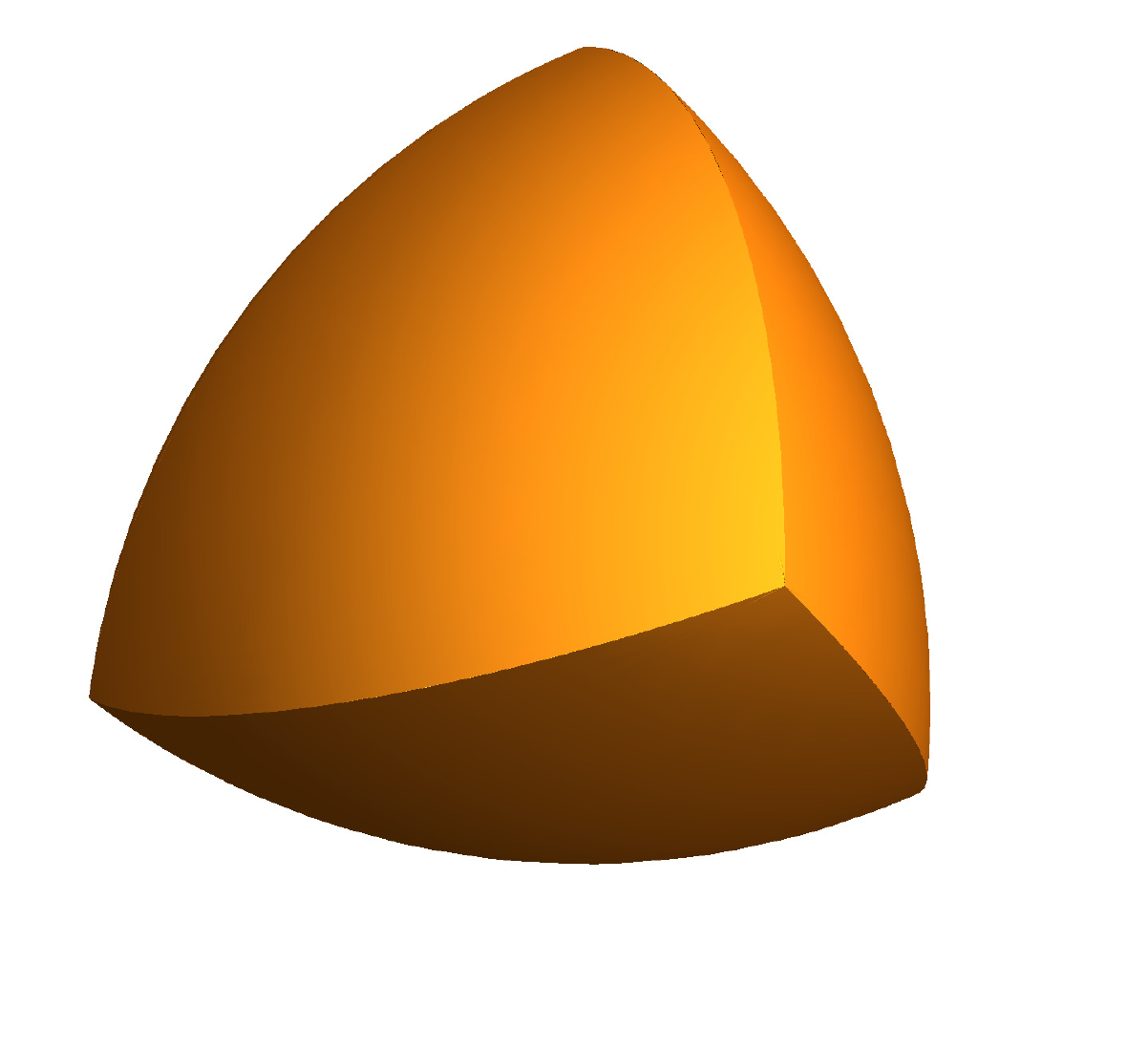}
     \end{subfigure}
          \caption{A Meissner tetrahedron having three rounded edges that form a triangle.}
        \label{SecondMeissnerFigure}
\end{figure}

\subsection{Notation and preliminaries}\label{NotAndPreSubSection}
We will use the term ``ball" to mean a closed ball and write
$$
B_r(a)=\{x\in \R^3: |x-a|\le r\}.
$$
\par {\bf The space $C(\Stwo)/\Pee$}.  We will denote $C(\Stwo)$ as the space of continuous $g: \Stwo\rightarrow \R$ endowed with the supremum norm
$$
\|g\|=\max_{|u|=1}|g(u)|.
$$ 
Let us also write $\Pee\subset C(\Stwo)$ for the closed subspace of functions 
$$
\Stwo\rightarrow \R; u\mapsto a\cdot u
$$
where $a$ ranges over points in $\R^3$. This corresponds to convex bodies which are singletons.  Since the volume and circumradius are invariant under translations, we will consider the quotient space 
$$
C(\Stwo)/\Pee=\{g+ \Pee: g\in C(\Stwo)\},
$$
which is endowed with the quotient norm 
$$
\|g+\Pee\|=\inf_{a\in \R^3}\max_{|u|=1}|g(u)+a\cdot u|.
$$
For ease of notation, we will use the variable $\bg$ to denote an element $C(\Stwo)/\Pee$ whenever there is $g\in \bg$ so that
$$
\bg=g+\Pee.
$$

\par {\bf The space $\Pee^\perp$}. We recall that the continuous dual space of $C(\Stwo)/{\Pee}$ when endowed with operator norm is isometrically isomorphic to  
$$
\Pee^\perp=\{\xi\in M(\Stwo): \xi|_{\Pee}=0\}.
$$ 
Here $M(\Stwo)=C(\Stwo)^*$ is the space of signed Radon measures on $\Stwo$ endowed with the total variation norm. In particular, $\Pee^\perp$ admits the norm 
$$
\|\xi\|_{*}=\sup\left\{\langle \xi,g\rangle: \bg\in C(\Stwo)/{\Pee},\; \|\bg\|\le 1\right\},
$$
where  
\be\label{Dualpairing}
\langle \xi,g\rangle:=\int_{\Stwo}g d\xi
\ee
is the natural pairing between $\xi\in M(\Stwo)$ and $g\in C(\Stwo)$.

\par We will denote $\Jay$ as the subdifferential of 
$$
C(\Stwo)/\Pee\ni \bg\mapsto \frac{1}{2}\|\bg\|^2
$$
and $\Jay^*$ for the subdifferential of 
$$
\Pee^\perp\ni \xi\mapsto \frac{1}{2}\|\xi\|_*^2.
$$
By convex duality, $\xi \in \Jay(\bg)$ if and only if $\bg\in \Jay^*(\xi)$ if and only if 
\be\label{DualityFormula}
\|\bg\|^2=\langle \xi,g\rangle=\|\xi\|_*^2.
\ee
\par {\bf The space $\Cee$}. Another natural space for us to study is
\be\label{theSpaceCee}
\Cee:=\left\{g\in C(\Stwo): g+\frac{1}{2}=H|_{\Stwo}, \text{$H$ is the support function of a constant width body in $\R^3$}\right\}.
\ee
We will identify a useful compactness property of $\Cee$ below.  And as mentioned above, if $K$ is a constant width body associated with $g\in\Cee$, then its circumradius is given by
$$
R(K)=\frac{1}{2}+\|\bg\|.
$$

\par In addition, we will consider the functional 
defined as
\par \be
E(\bg):=
\begin{cases}
\displaystyle\int_{\Stwo}\left(\frac{1}{2}|\nabla g|^2-g^2 \right)d\sigma, \quad &g\in\Cee\\
+\infty, \quad & g\not\in\Cee
\end{cases}
\ee
for $\bg\in C(\Stwo)/\Pee$.  We note that $E$ is well-defined since the integral $\int_{\Stwo}\left(\frac{1}{2}|\nabla g|^2-g^2 \right)d\sigma$ is invariant under the translation of $g$ by elements in $\Pee$. Moreover, 
it is not hard to see that $E$ is convex, proper, and lower-semicontinuous.   And as we previously noted,
$$
V(K)=\frac{\pi}{6}-\frac{1}{2}E(\bg)
$$
for a constant width body $K$ associated with $\bg$.  Since $V(K)\ge 0$, we also have
\be\label{BasicBoundonE}
\sup_{\bg\in \Cee}E(\bg)\le \frac{\pi}{3}.
\ee

\par Notice that if $g_1,g_2\in \Cee$ with corresponding constant with bodies $K_1,K_2\subset \R^3$,  then
$$
\|\bg_1-\bg_2\|=\inf_{a\in \R^3}d_{\mathcal H}(K_1+a,K_2).
$$
Here $d_{\mathcal H}$ is the Hausdorff distance.  In particular, if $\bg_n\rightarrow \bg$ and if $K_n$ and $K$ are respective convex bodies associated with $g_n\in\Cee$ and $g\in \Cee$, there is a sequence $(a_n)_{n\in \N}\subset\R^3$ for which 
$$
\lim_{n\rightarrow\infty}\|\bg_n-\bg\|=\lim_{n\rightarrow\infty}d_{\mathcal H}(K_n+a_n,K)=0.
$$
That is, up to translations, $K_n$ converges to $K$ in the Hausdorff topology.

\subsection{A doubly nonlinear evolution}
In what follows, we will study solutions $\xi: [0,\infty)\rightarrow\Pee^\perp$ of the {\it doubly nonlinear evolution}
\be\label{DNExi}
\partial E^*(\dot\xi(t))+\Jay^*(\xi(t))\ni 0\quad \text{a.e. $t\ge 0$}
\ee 
for a given initial condition.  Here 
$$
E^*(\zeta):=\sup\{\langle \zeta,g\rangle -E(\bg): g\in \Cee\}, \quad (\zeta\in \Pee^\perp)
$$
is the convex dual of $E$ and we recall  $\Jay^*$ is the subdifferential of $\frac{1}{2}\|\cdot\|_*^2$. It is typical to consider a solution of \eqref{DNExi} to mean that $\xi$ is absolutely continuous and that there is a measurable $\bg: [0,\infty)\rightarrow C(\Stwo)/\Pee$ 
with 
\be\label{geeacheconditions}
\bg(t)\in \Jay^*(\xi(t))\cap (-\partial E^*(\dot\xi(t)))
\ee
for almost every $t\ge 0$.   Let us for the moment suppose we have such a solution.
\\
\par {\bf First monotonicity formula}. By direct computation, we find
$$
\frac{d}{dt}\frac{1}{2}\|\xi(t)\|_*^2=-\left[E^*(\dot\xi(t))+E(\bg(t))\right]
$$
for almost every $t\ge 0$. Integrating this formula on the interval $[s,t]$ gives 
\be
\frac{1}{2}\|\xi(s)\|_*^2=\frac{1}{2}\|\xi(t)\|_*^2+\int^t_sE^*(\dot \xi(\tau))d\tau+\int^t_sE(\bg(\tau))d\tau.
\ee
This identity, combined with \eqref{DualityFormula} and \eqref{geeacheconditions}, implies
$$
\|\bg(t)\|=\|\xi(t)\|_*\text{  is a nondecreasing function of t.}
$$
\\
\par{\bf Second monotonicity formula}. In view of \eqref{geeacheconditions}, we also have $-\dot\xi(t)\in \partial E(\bg(t))$ for almost every $t\ge 0$. If, in addition $g:[0,\infty)\rightarrow C(\Stwo); t\mapsto g(t)$ is differentiable almost everywhere, then
\begin{align*}
\frac{d}{dt}E(\bg(t))&=-\langle  \dot \xi(t), \dot g(t)\rangle
\end{align*}
for almost every $t\ge 0$.  Using \eqref{DualityFormula}, we find
\begin{align}
&\left\langle \xi(t+\tau)-\xi(t), g(t+\tau)-g(t)\right\rangle \\
&=\|\xi(t+\tau)\|_*^2-\langle  \xi(t),g(t+\tau)\rangle-\langle  \xi(t+\tau),g(t)\rangle+\|\xi(t)\|_*^2 \\ 
&\ge  \|\xi(t+\tau)\|_*^2- \|\xi(t)\|_*\|\bg(t+\tau)\|-\|\xi(t+\tau)\|_*\|\bg(t)\|+\|\xi(t)\|_*^2 \\ 
&=  \|\xi(t+\tau)\|_*^2- 2\|\xi(t)\|_*\|\xi(t+\tau)\|_*+\|\xi(t)\|_*^2 \\ 
&=\left(\|\xi(t+\tau)\|_*-\|\xi(t)\|_*\right)^2.
\end{align}
Thus, we expect 
\be\label{SecondApriori}
\frac{d}{dt}E(\bg(t))=-\langle  \dot \xi(t), \dot g(t)\rangle\le -\left(\frac{d}{dt}\|\xi(t)\|_*\right)^2
\ee
for almost every $t\ge 0$. In particular,
$$
E(\bg(t)) \text{  is a nonincreasing function of $t$.}
$$
\\
\par {\bf Large time limits}. Since $E(\bg(t))$ is nonnegative, nondecreasing, and integrable on $[0,\infty)$, it must be that 
\be\label{preVolLim}
\lim_{t\rightarrow\infty}tE(\bg(t))=0.
\ee
Using this limit and the compactness of $\Cee$, we would then be able to conclude  
\be\label{preCircumRadLim}
\lim_{t\rightarrow\infty}\|\bg(t)\|=0.
\ee
\\
\par {\bf Geometric interpretation}. For each $t\ge 0$, $g(t)+1/2$ is the $\Stwo$ restriction of the support function of a constant width body $K_t\subset \R^3$. Given the monotonicity formulae above and the way circumradius and volume can be expressed in terms of the support function, 
$$
\text{$R(K_t)$ is a nonincreasing function of $t$}
$$
and 
$$
\text{$V(K_t)$ is a nondecreasing function of $t$}.
$$
In view of the above large time limits, 
$$
\begin{cases}
\displaystyle\lim_{t\rightarrow\infty}t\left(\frac{\pi}{6}-V(K_t)\right)=0\\\\
\displaystyle\lim_{t\rightarrow\infty}R(K_t)=\frac{1}{2}.
\end{cases}
$$
And up to translations, 
$$
\text{$K_t$ converges to $B_{1/2}(0)$ in the Hausdorff topology.} 
$$

\par {\bf Approach to existence}. Our goal is to establish that, for given initial conditions $g^0\in \Cee$ and $\xi^0\in \Pee^\perp$ such that
$$
\xi^0\in \Jay(\bg^0),
$$
there is a solution of  \eqref{DNExi} $\xi:[0,\infty)\rightarrow\Pee^\perp$ as described above.  Then we could 
obtain a mapping $\bg:[0,\infty)\rightarrow C(\Stwo)/\Pee$ and make the analogous geometric conclusions.  
However, we do not know how to carry this out primarily because the functional $E^*$ is not coercive. In particular, having an a priori bound on the integral
\be\label{EstarL1norm}
\int^\infty_0E^*(\dot\xi(t))dt
\ee
does not suggest that we can construct an absolutely continuous solution $\xi$. 

\par This is a typical problem encountered in the study of ``rate-independent" doubly nonlinear evolutions \cite{MR4158534,MR4197283,MR2210284,MR2105969,MR2290410,MR2887927,MR3531671,MR3740380,MR3636535}. For these flows, $E^*$ is usually a norm, so the gradient or subdifferential of $E^*$ is homogeneous of degree $0$.  It turns out that solutions to these flows have bounded variation.   While our functional $E^*$ is not a norm, it is convex, lower-semicontinuous and proper.  Using this basic information, we will develop the notion of the $E^*$ variation of a mapping $\xi:[0,\infty)\rightarrow\Pee^\perp$. Replacing $E^*$ with a norm in this definition results in the usual variation of $\xi$; and if $\xi$ is absolutely continuous, the $E^*$ variation of $\xi$ on $[0,\infty)$ is equal to the integral \eqref{EstarL1norm}.

\par These considerations lead to a notion of weak solution based on the ``classical" notion of solution described above. We will show that classical solutions are weak solutions. Moreover, we will prove that weak solutions $\xi$ exist for given initial conditions and their companion mappings $\bg$ essentially satisfy the properties discussed above.  To this end, we will use a compactness based approach by designing an approximation sequence,  establishing various bounds on these approximations, and then extracting a subsequence which converges to a solution with desired properties.  In order to carry out this procedure, we will need to first develop some ideas for support functions and to discuss the appropriate function spaces. These are the topics of the subsequent sections of this paper.

\section{Support functions}\label{SuppurtFunctionSect}
Let us recall a few basic facts about convex bodies.  Suppose $K\subset \R^3$ is a convex body with support function
\be
H(u)=\max_{x\in K}x\cdot u\quad (u\in \R^3).
\ee
Observe that $H$ is positively homogeneous and convex. Moreover, 
\be\label{setKinTermsofH}
K=\bigcap_{|u|=1}\left\{x\in \R^3: x\cdot u\le H(u)\right\}.
\ee
This identity can be used to show that the class of positively homogeneous, convex functions $H: \R^3\rightarrow \R$ are in one-to-one 
correspondence with convex bodies. 

\par If $K_1$ and $K_2$ are two convex bodies with respective 
support functions $H_1$ and $H_2$, we may consider their Minkowski sum 
$$
K_1+K_2:=\{x_1+x_2\in \R^3: x_1\in K_1,x_2\in K_2\},
$$
which is also convex body with support function $H_1+H_2$. In addition, we note that $K_1\subset K_2$ if and only if $H_1\le H_2$. 

\par The Hausdorff distance 
$$
d_{\mathcal H}(K_1,K_1):=\inf\{r>0: K_1\subset K_2+B_r(0),\; K_2\subset K_1+B_r(0) \}
$$
between $K_1$ and $K_2$ may also be expressed in terms of their support functions. Indeed 
\begin{align}
d_{\mathcal H}(K_1,K_1)&=\inf\left\{r>0: H_{K_1}(u)\le H_{K_2}(u)+r \text{ and }H_{K_1}(u)\le H_{K_2}(u)+r \text{ for $|u|=1$}\right\}\\
&=\inf\left\{r>0: \max_{|u|=1}|H_{K_1}(u)-H_{K_2}(u)|\le r \right\}\\
&=\max_{|u|=1}|H_{K_1}(u)-H_{K_2}(u)|.
\end{align}
Many more properties of support functions  can be found in the standard references on convex bodies such as \cite{MR920366, MR1491362,MR1216521}.

\par In this section, we will study a fixed a constant width body $K$ with support function $H$. We note that the equality
\be\label{ConstantWidthSumFormula}
K+(-K)=B_1(0),
\ee
is equivalent to $K$ having constant width. This in turn holds if and only if 
\be\label{ConstantWidthH}
H(u)+H(-u)=|u|
\ee
for each $u\in \R^3$.  Below, we will derive a formula for the circumradius and volume of $K$ in terms of $H$. We will 
also derive a few estimates on $H$ which will imply a compactness property for constant width bodies.  

\par We note that most if not all of the following results are likely to be found in the literature on constant width bodies (including \cite{MR920366,MR0123962,MR3930585}). However, we have included them in attempt to present a somewhat unified treatment of the support functions of constant width bodies.  Moreover, many of the results are valid in arbitrary dimension. Nevertheless, we will only focus on constant width bodies in $\R^3$. 

\subsection{In and circumradius}
We define the {\it inradius} of $K$ to be
$$
r(K)=\sup\{r\ge 0: K\supset B_r(a)\text{ some $a\in K$}\},
$$
and an {\it inball} as any ball with $B_{r(K)}(a)\subset K$.  Namely, an inball is a ball of maximal radius which can be included in $K$. 

\begin{prop}
$$
r(K)=\max_{a\in \R^3}\min_{|u|=1}\{H(u)-a\cdot u\},
$$
and the maximum occurs at some $a\in K$.  Moreover, $B_{r(K)}(a)$ is an inball.
\end{prop}
\begin{proof}
Note that $u\mapsto a\cdot u+r|u|$ is the support function of $B_r(a)$. Therefore,
\begin{align}
r(K)&=\sup\{r\ge0: K\supset B_r(a)\text{ some $a\in K$}\}\\
&=\sup\{r\ge0: a\cdot u+r\le H(u)\text{ all $|u|=1$ and some $a\in K$}\}\\
&=\sup\left\{r\ge0: r\le\min_{|u|=1}\{H(u)-a\cdot u\}\text{ for some $a\in K$}\right\}\\
&=\sup\left\{r\ge0: r\le\max_{a\in K}\min_{|u|=1}\{H(u)-a\cdot u\}\right\}\\
&=\max_{a\in K}\min_{|u|=1}\{H(u)-a\cdot u\}.
\end{align}
If we choose $a\in K$ so that  
$$
r(K)=\min_{|u|=1}\{H(u)-a\cdot u\},
$$
then $r(K)+a\cdot u\le H(u)$ for all $|u|=1$. That is, $B_{r(K)}(a)\subset K$.

\par We next claim that, $r(K)$ is equal to  
$$
\overline{r}(K):=\sup_{a\in \R^3}\min_{|u|=1}\{H(u)-a\cdot u\}.
$$
So far we have $r(K)\le \overline{r}(K)$.  Moreover, for any $a\in \R^3$ there is a unit vector $u$ for which $a\cdot u\ge 0$. Thus
$$
\min_{|u|=1}\{H(u)-a\cdot u\}\le \max_{|u|=1}H(u)<\infty.
$$
It follows that $\overline{r}(K)<\infty$. We may suppose $ \overline{r}(K)>0$ or else  $ \overline{r}(K)=r(K)=0$. In this case, let $\epsilon\in (0,\overline{r}(K))$ and choose $a^\epsilon$ such that 
$$
\overline{r}(K)-\epsilon\le \min_{|u|=1}\{H(u)-a^\epsilon\cdot u\}.
$$
Then $B_{\overline{r}(K)-\epsilon}(a^\epsilon)\le K$. In particular, $a^\epsilon\in K$ and 
$$
\overline{r}(K)-\epsilon\le r(K).
$$
\end{proof}

\par Analogously, we can define the {\it circumradius} of $K$ as
$$
R(K)=\inf\{r>0: K\subset B_r(a)\text{ some $a\in \R^3$}\}.
$$
 Note that since $K$ is compact, $R(K)<\infty$. We'll also call a {\it circumball} a ball such that $K\subset B_{R(K)}(a)$. 
 First we show that circumballs are unique. 

\begin{lem}
There can be only one $a\in \R^3$ such that $K\subset B_{R(K)}(a)$.
\end{lem}
\begin{proof}
Set $r:=r(K)$ and suppose $K\subset B_{r}(a_1)\cap B_{r}(a_2)$ with $a_1\neq a_2$. For any $z\in K$, $|z-a_1|\le r$ and $|z-a_2|\le r$. Therefore, 
$$
|a_1-a_2|< |a_1-z|+|z-a_2|\le 2r. 
$$
Set 
$$
s:=\sqrt{r^2-\left|\frac{a_1-a_2}{2}\right|^2},
$$
and note that $0< s<r$.  Moreover, for $z\in K$
$$
\left|z-\frac{a_1+a_2}{2}\right|^2=\frac{1}{2}|z-a_1|^2+\frac{1}{2}|z-a_2|^2-\left|\frac{a_1-a_2}{2}\right|^2\le r^2-\left|\frac{a_1-a_2}{2}\right|^2=s^2.
$$
Thus, $K\subset B_{s}((a_1+a_2)/2).$ This contradicts the assumption that $r=r(K)$. Thus, it must be that $a_1=a_2$.

\end{proof}
In analogy with our formula for the inradius, we have the following formula for the circumradius.  
\begin{prop} 
\be\label{RCircumFormula}
R(K)=\min_{a\in \R^3}\max_{|u|=1}\{H(u)-a\cdot u\}.
\ee
Moreover, $B_{R(K)}(a)$ is the circumball provided that
$$
R(K)=\max_{|u|=1}\{H(u)-a\cdot u\}.
$$ 
\end{prop}
\begin{proof}
Observe that if $K\subset B_r(a)$ for  some $a\in \R^3$, then $H(u)\le a\cdot u+r$ for all $|u|=1$. As a result, 
\begin{align}
r&\ge \max_{|u|=1}\{H(u)-a\cdot u\}\\
 &\ge \inf_{a\in \R^3}\max_{|u|=1}\{H(u)-a\cdot u\}\\
 &=: \overline{R}(K).
\end{align}
Thus $R(K)\ge  \overline{R}(K).$ Now let $\epsilon>0$ and select $a^\epsilon\in \R^3$ for which 
$$
\overline{R}(K)+\epsilon\ge \max_{|u|=1}\{H(u)-a^\epsilon\cdot u\}.
$$
This implies, $K\subset B_{\overline{R}(K)+\epsilon}(a^\epsilon)$. Thus, 
$$
R(K)\le \overline{R}(K)+\epsilon.
$$
It follows that $R(K)= \overline{R}(K).$

\par Suppose $(a^k)_{k\in \N}$ is a sequence in $\R^3$ for which 
$$
R(K)=\lim_{k\rightarrow\infty}\max_{|u|=1}\{H(u)-a^k\cdot u\}.
$$
Then for all $|v|=1$ and sufficiently large $k\in \N$, 
\begin{align}
-a^k\cdot v&=H(v)-a^k\cdot v -H(v)\\
&\le \max_{|u|=1}\{H(u)-a^k\cdot u\}+\max_{|u|=1}|H(u)|\\
&\le R(K)+1+\max_{|u|=1}|H(u)|.
\end{align}
This implies that $(a^k)_{k\in \N}$ is bounded and therefore has a convergent subsequence $(a^{k_j})_{j\in \N}$ with limit $a^\infty$. As a result, 
$$
R(K)=\lim_{j\rightarrow\infty}\max_{|u|=1}\{H(u)-a^{k_j}\cdot u\}
=\max_{|u|=1}\{H(u)-a^\infty \cdot u\}.
$$
We conclude \eqref{RCircumFormula}. Since $H(u)\le R(K)+a^\infty \cdot u$ it must also be that $K\subset B_{R(K)}(a^\infty)$.
\end{proof}

\par We will now make use of \eqref{ConstantWidthH}.  
\begin{lem}
$K$ has a unique inball $B_{r(K)}(a)$ which is concentric with its circumball.  Moreover, 
$$
r(K)+R(K)=1.
$$
\end{lem}
\begin{proof}
Since $H(u)+H(-u)=1$ for each $|u|=1$, 
\begin{align}
r(K)&=\max_{a\in \R^3}\min_{|u|=1}\{H(u)-a\cdot u\}\\
&=\max_{a\in \R^3}\min_{|u|=1}\{1-H(-u)-a\cdot u\}\\
&=1+\max_{a\in \R^3}\min_{|u|=1}\{-H(-u)-a\cdot u\}\\
&=1-\min_{a\in \R^3}\max_{|u|=1}\{H(-u)+a\cdot u\}\\
&=1-\min_{a\in \R^3}\max_{|u|=1}\{H(u)-a\cdot u\}\\
&=1-R(K).
\end{align}
Now select $a\in K$ for which
$$
r(K)=\min_{|u|=1}\{H(u)-a\cdot u\}.
$$
This implies $H(u)\ge a\cdot u+r(K)$ for $|u|=1$ so that 
 $B_{r(K)}(a)\subset K$.  In addition, we have 
\begin{align}
R(K)&=1-r(K)\\
&=1-\min_{|u|=1}\{H(u)-a\cdot u\}\\
&=\max_{|u|=1}\{1-H(u)+a\cdot u\}\\
&=\max_{|u|=1}\{H(-u)+a\cdot u\}\\
&=\max_{|u|=1}\{H(u)-a\cdot u\}.
\end{align}
Thus, $R(K)+a\cdot u\ge H(u)$ for all $|u|=1$. Furthermore, $K\subset B_{R(K)}(a)$ and $a$ is uniquely specified. 
\end{proof}
\begin{cor}\label{UpperBoundRKay}
$$
R(K)\le\sqrt{\frac{3}{8}}\quad \text{and}\quad r(K)\ge 1-\sqrt{\frac{3}{8}}.
$$
Equality holds in either inequality if $K$ is a Meissner tetrahedron. 
\end{cor}
\begin{proof}
In view of \eqref{ConstantWidthSumFormula}, $K$ has diameter 1. Jung's theorem implies
$$
R(K)\le \sqrt{\frac{3}{2(3+1)}}=\sqrt{\frac{3}{8}}
$$ 
and that equality would hold if $K$ is a regular tetrahedron.  Now let $M\subset \R^3$ be a Meissner tetrahedron. Then $M$ includes a regular tetrahedron $T$ of diameter 1. Therefore, 
$$
R(M)\ge R(T)=\sqrt{\frac{3}{8}}.
$$
By the previous lemma, we also have
$$
r(K)=1-R(K)\ge 1-\sqrt{\frac{3}{8}}
$$
and similarly conclude that equality holds if $K$ is a Meissner tetrahedron.
\end{proof}
We can now derive a key formula for the circumradius of a body of constant width.
\begin{prop}\label{propFormulafofR}
$$
R(K)=1/2+\min_{a\in \R^3}\max_{|u|=1}|H(u)-1/2-a\cdot u|
$$
and the minimum occurs at the center of the circumball for $K$.
\end{prop}
\begin{proof}
Let $a\in K$ be the center of $K$'s circumball. Then 
$$
a\cdot u+r(K)\le H(u)\le a\cdot u+R(K)
$$
for all $|u|=1$. As $r(K)+R(K)=1$, 
$$
|H(u)-1/2-a\cdot u|\le R(K)-1/2.
$$
Thus, 
$$
\max_{|u|=1}|H(u)-1/2-a\cdot u|\le R(K)-1/2.
$$
Also note 
\begin{align}
R(K)-1/2&=\inf_{b\in \R^3}\max_{|u|=1}\{H(u)-1/2-b\cdot u\}\\
&\le \inf_{b\in \R^3}\max_{|u|=1}|H(u)-1/2-b\cdot u|.
\end{align}
\end{proof}
\begin{rem}
A corollary of Proposition \ref{propFormulafofR} is 
\be\label{propFormulafofR2}
R(K)=\frac{1}{2}+d_{{\mathcal H}}(K-a,B_{1/2}(0)),
\ee
where $a$ is $K$'s circumball.
\end{rem}
\begin{cor}\label{EstimatesSuppHnoted}
If $B_{R(K)}(a)$ is the circumball of $K$, then
\begin{enumerate}[(i)]

\item  $|H(u)-1/2-a\cdot u|\le\displaystyle \sqrt{\frac{3}{8}}-\frac{1}{2}$  for $|u|=1$, and

\item $|H(u)-H(v)-a\cdot (u-v)|\le |u-v|$  for $u,v\in \R^3$.

\end{enumerate}
\end{cor}
\begin{proof}
$(i)$ This follows from Corollary \ref{UpperBoundRKay} and Proposition \ref{propFormulafofR}. $(ii)$ Choose $x\in K$ for which $H(u)=x\cdot u$.  As the diameter of $K$ is 1,  
\begin{align}
H(u)-H(v)-a\cdot (u-v)&\le x\cdot u-x\cdot v-a\cdot (u-v)\\ 
&=(x-a)\cdot (u-v)\\
&\le |u-v|.
\end{align}
The same argument implies $H(v)-H(u)-a\cdot (v-u)\le  |u-v|.$ 
\end{proof}
\subsection{A gradient estimate}
Since $H$ is convex, it is twice differentiable almost everywhere by Rademacher's theorem.  In view of \eqref{ConstantWidthH}, we additionally have 
\be\label{W2infinity}
0\le D^2H(u)w\cdot w\le \frac{1}{|u|}
\ee
for almost every $u\in \R^3\setminus\{0\}$ and each $w\in \Stwo$. This in turn implies that $H$ is continuously differentiable away from $0$ \cite{MR2233133}.  We can also quantify this below with the following estimate.
\begin{prop}
For $|u|,|v|\ge r$,
$$
|DH(u)-DH(v)|\le \displaystyle\frac{1}{r}\left(1+\frac{\pi}{2}\right)|u-v|.
$$
\end{prop}
\begin{proof}
We will first prove the estimate assuming $|u|=|v|=1$, then $|u|=|v|=r$, and finally $|u|,|v|\ge r$. To this end,  let $H^\epsilon=\eta^\epsilon*H$ where $\eta^\epsilon$ is a standard mollifier. That is, $\eta\in C_c^\infty$ is nonnegative, the support of $\eta$ is the unit ball in $\R^3$, with $\int_{\R^3}\eta(y)dy=1$, 
and $\eta^\epsilon(x)=\eta(x/\epsilon)/\epsilon^3$ for $\epsilon>0$. It is not hard to see that $H^\epsilon$ is smooth and also converges in $C^1_{\text{loc}}(\R^3\setminus\{0\})$ to $H$ (see Chapter 4 of \cite{MR3409135} for example).  Let $u,v\in\Stwo$ and choose a geodesic path $\gamma: [0,1]\rightarrow\Stwo$ joining $u$ to $v$.  Then 
\be
DH^\epsilon(v)-DH^\epsilon(u)=\int^1_0D^2H^\epsilon(\gamma(t))\dot\gamma(t)dt.
\ee
Observe that for $\epsilon<1$ and $w\in \Stwo$, 
\begin{align}
D^2H^\epsilon(u)w\cdot w&=\int_{B_\epsilon(0)}\eta^\epsilon(x)D^2H(u-x)w\cdot wdx\\
&\le \int_{B_\epsilon(0)}\eta^\epsilon(x)\frac{1}{|u-x|}dx\\
&\le \int_{B_\epsilon(0)}\eta^\epsilon(x)\frac{1}{|u|-|x|}dx\\
&\le \int_{B_\epsilon(0)}\eta^\epsilon(x)\frac{1}{1-\epsilon}dx\\
&=\frac{1}{1-\epsilon}.
\end{align}
\par It follows that
\begin{align}
|DH^\epsilon(v)-DH^\epsilon(u)|&\le \int^1_0|D^2H^\epsilon(\gamma(t))||\dot\gamma(t)|dt\le \frac{1}{1-\epsilon}\int^1_0|\dot\gamma(t)|dt= d(u,v) \frac{1}{1-\epsilon}
\end{align}
where 
$$
d(u,v):=\cos^{-1}(u\cdot v)
$$
is the standard metric on $\Stwo$. As
\be
d(u,v)\le \frac{\pi}{2}|u-v|
\ee
for $u,v\in \Stwo$, we may send $\epsilon\rightarrow0$ and arrive at 
\begin{align}
|DH(u)-DH(v)|\le  \frac{\pi}{2}|u-v|.
\end{align}
When $|u|=|v|=r$, we use the fact that $DH$ is $0$-degree homogeneous on $\R^3\setminus\{0\}$ to get 
$$
|DH(u)-DH(v)|=|DH(u/r)-DH(v/r)|\le  \frac{\pi}{2}|u/r-v/r|= \frac{\pi}{2r}|u-v|.
$$

\par Now let $|v|\ge |u|\ge r$. By the triangle inequality 
\be\label{TriangleInequalityDH}
|DH(v)-DH(v)|\le \left|DH(v)-DH\left(|u|\frac{v}{|v|}\right)\right|+\left|DH\left(|u|\frac{v}{|v|}\right)-DH(u)\right|.
\ee
We can estimate the second term from the computation above: 
$$
\left|DH\left(|u|\frac{v}{|v|}\right)-DH(u)\right|\le \frac{\pi}{2r}\left|u-|u|\frac{v}{|v|}\right|.
$$
An elementary computation shows 
$$
\left|u-|u|\frac{v}{|v|}\right|\le |u-v|;
$$
in general, this inequality holds provided $|v|\ge |u|$ with $|v|>0$. Therefore, 
$$
\left|DH\left(|u|\frac{v}{|v|}\right)-DH(u)\right|\le \frac{\pi}{2r}|u-v|.
$$
\par As for the first term in \eqref{TriangleInequalityDH}, we employ the linear path
$$
\alpha(t)=|u|\frac{v}{|v|}+t\left(v-|u|\frac{v}{|v|}\right), \quad 0\le t\le 1.
$$
Note in particular that 
$$
|\alpha(t)|\ge\alpha(t)\cdot \frac{v}{|v|}=|u|+t(|v|-|u|)\ge |u|\ge r.
$$
Consequently, if $H$ is smooth
\begin{align}
\left|DH(v)-DH\left(|u|\frac{v}{|v|}\right)\right|&=\left|DH(\alpha(1))-DH\left(\alpha(0)\right)\right|\\
&\le \int^1_0|D^2H(\alpha(t))||\dot\alpha(t)|dt\\
&\le \left|v-|u|\frac{v}{|v|}\right|\int^1_0\frac{1}{|\alpha(t)|}dt\\
&\le \frac{|v|-|u|}{r}\\
&\le \frac{|u-v|}{r}.
\end{align}
Otherwise, we can smooth $H$ with a mollifier and derive the same inequality. 
Putting the two bounds together gives 
$$
|DH(u)-DH(v)|\le \displaystyle\frac{1}{r}\left(1+\frac{\pi}{2}\right)|u-v|
$$
for $|u|,|v|\ge r$.
\end{proof}
\begin{cor}\label{CompactnessHemm}
Suppose for each $m\in \N$ that $H^m$ is the support function of constant width body $K^m\subset \R^3$ with circumball centered at $a^m$. There is a subsequence 
$(H^{m_j})_{j\in \N}$ and a support function $H^\infty$ of a constant width body $K^\infty\subset \R^3$ such that 
\be\label{HemJayUNif}
\lim_{j\rightarrow\infty}\max_{|u|\le s}|H^{m_j}(u)-a^{m_j}\cdot u-H^\infty(u)|=0
\ee
and
\be\label{DHemJayUNif}
\lim_{j\rightarrow\infty}\max_{r\le |u|\le 1/r}|DH^{m_j}(u)-a^{m_j}-DH^\infty(u)|=0
\ee
for each $s>0$ and $0<r\le 1$.
\end{cor}
\begin{proof}
Set
$$
\tilde H^m(u):=H^{m}(u)-a^{m_j}\cdot u
$$
for $u\in \R^3$.  By Corollary \ref{EstimatesSuppHnoted}, 
$$
|\tilde H^m(u)-\tilde H^m(v)|\le|u-v|
$$
for all $u,v\in \R^3$ and each $m\in \N$. Since $H^m(0)=0$ for each $m\in \N$,  the Arzel\`a-Ascoli theorem implies $\tilde H^{m_j}$ has a locally uniformly convergent subsequence.  The limit function $H^\infty$ is necessarily continuous, sublinear, and satisfies the constant width condition \eqref{ConstantWidthH}. As a result, it is the 
support function of the constant width body 
$$
K^\infty:=\bigcap_{|u|=1}\left\{x\in \R^3: x\cdot u\le H^\infty(u)\right\}.
$$
This establishes the first limit \eqref{HemJayUNif}. As for the second limit, we recall 
$$
| D\tilde H^m(u)- D\tilde H^m(v)|\le \frac{1}{r}\left(1+\frac{\pi}{2}\right)|u-v|
$$
for $r\le |u|,|v|\le 1/r$ for each $0<r\le 1$. It follows that $(\tilde H^{m_j})_{j\in \N}$ converges in $C^1_{\text{loc}}(\R^3\setminus\{0\})$ to $H^\infty$.  This establishes
 \eqref{DHemJayUNif}.
\end{proof}
\begin{rem}
The limit \eqref{HemJayUNif} implies that 
$$
K^{m_j}+a^{m_j}\rightarrow K^\infty
$$
in the Hausdorff topology. This could have also been deduced as a consequence of Blaschke's selection theorem (Theorem 2.5.2 in \cite{MR3930585}). Moreover, the argument above 
can be easily adapted to fashion an elementary proof of Blaschke's selection theorem. 
\end{rem}

\subsection{Formula for the volume}
Let us set
$$
h:=H|_{\Stwo}
$$
as the restriction of $H$ to the unit sphere in $\R^3$. Since $H$ is positively homogeneous, this function encodes all information about $H$. It also naturally 
inherits various properties from $H$ such as
$h\in C^{1,1}(\Stwo)$ and 
\be\label{LittleHConstWidth}
h(u)+h(-u)=1
\ee
for each $u\in\Stwo$.

\par We can also infer information on the gradient and Hessian of $h:\Stwo\rightarrow\R$.  The following assertion was proved by Howard along the way to his solution of Nakajima's problem in $\R^3$ \cite{MR2233133}. 

\begin{lem}\label{HowardLemma}
(i) For each $u\in\Stwo$,
\be
DH(u)=\nabla h(u)+h(u)u\in \partial K
\ee
and $u$ is an outward normal vector to $\partial K$ at $DH(u)$. \\
(ii) The map $DH: \Stwo\rightarrow \partial K$ is Lipschitz and surjective and its gradient is given by
\be\label{hHessianCondition}
\nabla(DH)(u)=\nabla^2h(u)+h(u)\id 
\ee
for $\sigma$ almost every $u\in \Stwo$. Here $u^\perp:=\{v\in \R^3: v\cdot u=0\}$. \\
(iii) The linear transformation
$$
\nabla(DH)(u):u^\perp\rightarrow u^\perp
$$
is nonnegative definite for $\sigma$ almost every $u\in \Stwo$. 
\end{lem}

\par Thus far, much of what we've done in $\R^3$ carries over with little to no change in $\R^n$. The following result however is specific to convex bodies in $\R^3$. 
The expression below for the volume of $K$ is known as {\it Blaschke's relation}.

\begin{prop} 
\be\label{BlasckeRelation}
V(K)=\frac{1}{2}\sigma(\partial K)-\frac{\pi}{3},
\ee
where
\be\label{VolumeFormula}
\sigma(\partial K)=\int_{\Stwo}\left(h^2-\frac{1}{2}|\nabla h|^2 \right)d\sigma.
\ee
\end{prop}
\begin{proof}
1. Let $t>0$ and consider $K^t:=K+B_t(0)$. This is a convex body with support function
$$
H^t(u):=H(u)+t|u|.
$$
Note that 
$$
DH^t(u)=DH(u)+tu=\nabla h(u)+h(u)u+tu
$$
for $u\in\Stwo$. It follows that
$$
(DH^t(u)-DH^t(v))\cdot (u-v)\ge t|u-v|^2
$$
and thus
$$
|DH^t(u)-DH^t(v)|\ge t|u-v|
$$
for $u,v\in \Stwo$. As a result, $DH^t: \Stwo\rightarrow \partial (K^t)$ is a bi-Lipschitz map. 

\par Let $u^t: \partial (K^t)\rightarrow \Stwo$ denote the inverse of $DH^t$; by part $(i)$ of Lemma \ref{HowardLemma} this is the outward unit normal field on $ \partial (K^t)$. Applying the divergence theorem, the area formula (Theorem 3.9 of \cite{MR3409135}), and parts $(ii)$ and $(iii)$ of Lemma \ref{HowardLemma}, we find 
\begin{align}
V(K^t)&=\int_{K^t}dx\\
&=\frac{1}{3}\int_{K^t}\text{div}(x)dx\\
&=\frac{1}{3}\int_{\partial (K^t)}x\cdot u^t(x)d\sigma(x)\\
&=\frac{1}{3}\int_{\Stwo }DH^t(u)\cdot u\;\text{det}(\nabla(DH^t)(u))d\sigma(u)\\
&=\frac{1}{3}\int_{\Stwo }H^t(u)\text{det}(\nabla(DH^t)(u))d\sigma(u)\\
&=\frac{1}{3}\int_{\Stwo }(h(u)+t)\text{det}(\nabla^2h(u)+(h(u)+t)\id)d\sigma(u).
\end{align}
Note that $K^s\subset K^t$ for $0< s\le t$ and $K=\bigcap_{t>0}K^t.$ We can then employ the monotonicity of Lebesgue measure and send $t\rightarrow 0^+$ to get
 $$
 V(K)=\frac{1}{3}\int_{\Stwo }h(u)\text{det}(\nabla^2h(u)+h(u)\id)d\sigma(u).
 $$
 
 \par 2. Similarly we find
 \be\label{sigmatformula}
 \sigma(\partial K^t)=\int_{\Stwo }\text{det}(\nabla^2h(u)+(h(u)+t)\id)d\sigma(u)
 \ee
for $t> 0$.  We claim this formula also holds at $t=0$. To see this, we let $r>0$ be the inradius of $K$ and observe  
$$
K=(1-t)K+tK\supset (1-t)K+t B_r(0)=(1-t)[K^{tr/(1-t)}]
$$
for $t\in (0,1)$. Since $\sigma(\partial A)\le \sigma(\partial B)$ for two convex bodies with $A\subset B$ \cite{MR3756927}, 
$$
\sigma(K)\ge (1-t)^2  \sigma(\partial K^{tr/(1-t)})= (1-t)^2\int_{\Stwo }\text{det}\left(\nabla^2h(u)+\left(h(u)+\frac{rt}{1-t}\right)\id\right)d\sigma(u)
$$
for $t\in (0,1$).  Likewise, $\sigma(K)$ is bounded from above by the right hand side of \eqref{sigmatformula} for each $t>0$. We conclude
$$
 \sigma(\partial K)=\int_{\Stwo }\text{det}(\nabla^2h(u)+h(u)\id)d\sigma(u).
 $$

\par We will now employ the coarea formula (Theorem 3.11 of \cite{MR3409135}) with the Lipschitz function
 $$
 f(x)=d(x,K).
 $$
 Note that $|Df(x)|=1$ almost everywhere in $K^c$. In particular, for $t\ge 0$
 \begin{align}
 V(K^t)- V(K)&=\int_{K^t\setminus K}dx\\
& =\int_{\R^3}1_{K^t\setminus K}(x)dx\\
& =\int_{\R^3}1_{K^t\setminus K}(x)|Df(x)|dx\\
& =\int^\infty_{0}\left\{\int_{f^{-1}(\{s\})}1_{K^t\setminus K}(x)d\sigma(x)\right\}ds\\
& =\int^\infty_{0}\sigma((K^t\setminus K)\cap f^{-1}(\{s\})ds\\
& =\int^t_{0}\sigma(\partial K^s)ds.
 \end{align}
 Here we used $d(x,K)=s>0$ if and only if $x\in \partial K^s$.   Thus, 
\begin{align}
\sigma(\partial K)&=\frac{d}{dt}V(K^t) \bigg|_{t=0}\\
&=\left.\frac{d}{dt}\frac{1}{3}\int_{\Stwo }(h(u)+t)\text{det}(\nabla^2h(u)+(h(u)+t)\id)d\sigma(u)\right|_{t=0}\\
&=\frac{1}{3}\int_{\Stwo }\text{det}(\nabla^2h(u)+h(u)\id)d\sigma(u)+\frac{1}{2}\int_{\Stwo}h(u)(\Delta h(u)+2h(u))d\sigma(u)\\
&=\frac{1}{3}\sigma(\partial K)+\frac{1}{3}\int_{\Stwo}h(u)(\Delta h(u)+2h(u))d\sigma(u).
\end{align}
That is, 
\begin{align}
\sigma(\partial K)&=\frac{1}{2}\int_{\Stwo}h(\Delta h+2h)d\sigma\\
&=\frac{1}{2}\int_{\Stwo}h\Delta h+2h^2d\sigma\\
&=\int_{\Stwo}\left(h^2-\frac{1}{2}|\nabla h|^2\right)d\sigma,
\end{align}
which is \eqref{VolumeFormula}. 

\par 3. As for Blaschke's relation, we will use the fact that Lebesgue measure is invariant under orthogonal transformations so that $V(K)=V(-K).$ 
Moreover, the support function of $-K$ is 
$$
H(-u)=|u|-H(u)\quad (u\in \R^3),
$$
as $K$ has constant width. Consequently, 
\begin{align}
&V(K)\\
&=V(-K)\\
&=\frac{1}{3}\int_{\Stwo }h(-u)\text{det}(\nabla^2h(-u)+h(-u)\id)d\sigma(u)\\
&=\frac{1}{3}\int_{\Stwo }(1-h(u))\text{det}(-\nabla^2h(u)+(1-h(u))\id)d\sigma(u)\\
&=\frac{1}{3}\int_{\Stwo }(1-h)\text{det}(1-(\nabla^2h+h\id))d\sigma\\
&=\frac{1}{3}\int_{\Stwo }(1-h)[1-(\Delta h+2h)+\text{det}(\nabla^2h+h\id)]d\sigma\\
&=\frac{1}{3}\int_{\Stwo }\left[1-(\Delta h+2h)+\text{det}(\nabla^2h+h\id) -h(1-(\Delta h+2h)+\text{det}(\nabla^2h+h\id)\right]d\sigma\\
&=\frac{1}{3}\int_{\Stwo }\left[1-(\Delta h+2h)+\text{det}(\nabla^2h+h\id) -h+h(\Delta h+2h)-h\text{det}(\nabla^2h+h\id)\right]d\sigma.
\end{align}

\par In view of the constant width condition \eqref{LittleHConstWidth}, 
$$
\frac{1}{4\pi}\int_{\Stwo}hd\sigma=\frac{1}{2}.
$$
Since $\Delta h=\text{div}(\nabla h)$ integrates to $0$ over $\Stwo$, it follows that 
$$
\int_{\Stwo }[1-(\Delta h+2h)]d\sigma=0.
$$
Using the formulae above that we derived for $V(K)$ and $\sigma(\partial K)$ then gives 
\begin{align}
&V(K)\\
&=\frac{1}{3}\int_{\Stwo }\left[1-(\Delta h+2h)+\text{det}(\nabla^2h+h\id) -h+h(\Delta h+2h)-h\text{det}(\nabla^2h+h\id)\right]d\sigma\\
&=\frac{1}{3}\sigma(\partial K)-\frac{2\pi}{3}+\frac{1}{3}\int_\Stwo h(\Delta h+2h)d\sigma -V(K)\\
&=\sigma(\partial K)-\frac{2\pi}{3} -V(K).
\end{align}
We conclude \eqref{BlasckeRelation}.
\end{proof}

\section{Spaces of functions and measures}\label{OurSpacesSect}
In this section, we will study the various spaces that will be needed in our analysis of the doubly nonlinear 
evolution \eqref{DNExi}. Some of these spaces were introduced in subsection \ref{NotAndPreSubSection}. 
 First, we will show that the space $\Cee$ is compact in a certain sense. Then 
we will consider a subspace of $H^1(\Stwo)$ with an inner product tailored for this 
work. Next, we will study measures in $\Pee^\perp$ modulo the ones which vanish on $\Cee$; we'll also identify an 
important compact subset of this quotient space that will play a crucial role in our subsequent existence proof. Finally, we will 
introduce the notion of the $E^*$ variation of a path $\xi: [0,\infty)\rightarrow\Pee^\perp$. 

\subsection{The space $\Cee$}
Let us recall the space $\Cee$ defined in \eqref{theSpaceCee}.  We will argue that this space is convex and inherits the compactness we discussed above for support functions.   

\begin{lem}\label{CompactCeeLem}
(i) If $g\in \Cee$, then 
$$
\|\bg\|\le \sqrt{\frac{3}{8}}-\frac{1}{2}.
$$
(ii) $\Cee$ is convex.\\
(iii) For any sequence $(g^k)_{k\in \N}\subset {\mathcal C}$, there is a subsequence $(\bg^{k_j})_{j\in \N}$ which converges to some $\bg$ with $g\in \Cee$. Moreover, there is $a^{j}\in \R^3$ for each $j\in \N$ such that $g^{k_j}(u)-a^j\cdot u$ converges in $C^1(\Stwo)$ to $g$.
\end{lem}
\begin{proof}
$(i)$ Suppose $g\in \Cee$. There is a constant width body such that $g+1/2$ is the $\Stwo$ restriction of the support function of $K\subset \R^3$. By Corollary \ref{UpperBoundRKay}, 
$$
R(K)=\frac{1}{2}+\|\bg\|\le \sqrt{\frac{3}{8}}.
$$

\par $(ii)$ Assume $g_0,g_1\in \Cee$ and let $H_0, H_1$ denote the corresponding support functions with 
respective constant width bodies $K_0,K_1\subset \R^3$. Suppose $\lambda\in [0,1]$ and note that 
the support function of the convex body $(1-\lambda) K_0+\lambda K_1$ is $(1-\lambda) H_0+\lambda H_1$. Since 
$$
(1-\lambda) H_0(u)+\lambda H_1(u)+[(1-\lambda) H_0(-u)+\lambda H_1(-u)]=|u|,
$$
$(1-\lambda) K_0+\lambda K_1$ has constant width. As 
$$
((1-\lambda) g_0+\lambda g_1)+1/2=(1-\lambda)(g_0+1/2)+\lambda (g_1+1/2)=[(1-\lambda) H_0+\lambda H_1]|_{\Stwo},
$$
we conclude $(1-\lambda)g_0+\lambda g_1\in \Cee$. Therefore, $\Cee$ is convex.

\par $(iii)$ Let  
$$
 g^k+1/2=H^k|_{\Stwo},
$$
where $H^k$ is the support function of a constant width body.  By Corollary \ref{CompactnessHemm}, there are $a^j\in \R^3$ and a subsequence  $(H^{k_j})_{j\in \N}$ for which $H^{k_j}(u)-a^j\cdot u$ converges uniformly to the support function of a constant width body $H$ on compact subsets of $\R^3$. If $g+1/2=H|_\Stwo$, 
$$
\|\bg^{k_j}-\bg\|\le \max_{|u|=1}| g^{k_j}(u)-g(u)-a^j\cdot u|=\max_{|u|=1}|H^{k_j}(u)-H(u)-a^j\cdot u|\rightarrow 0
$$
 as $j\rightarrow\infty$.
 
\par Corollary \ref{CompactnessHemm} also gives that $D H^{k_j}(u)-a^j$ converges to $DH(u)$ uniformly for $u\in \Stwo$.  This implies
$$
\nabla ( g^{k_j}(u)-a^j\cdot u)=DH^{k_j}(u)-a^j-( g^{k_j}(u)-a^j\cdot u+1/2)u\rightarrow DH(u)-(g(u)+1/2)u=\nabla g(u)
$$
uniformly for $u\in \Stwo$. 
\end{proof}

\subsection{Square integrable paths}
The Sobolev space $H^1(\Stwo)$ is the completion of $C^\infty(\Stwo)$ in the norm 
\be\label{standardNormH1}
g\mapsto \left(\int_\Stwo \left(g^2+|\nabla g|^2\right)d\sigma\right)^{1/2}.
\ee
We refer the reader to Chapter 2 of \cite{MR1481970} for more on $H^1(\Stwo)$.  We will consider the closed subspace
$$
{\mathcal V}:=\left\{g\in H^1(\Stwo): \int_{\Stwo}gd\sigma =0\right\}
$$
which admits the semi-inner product 
$$
(g,h):=\int_{\Stwo}(\nabla g\cdot \nabla h-2gh)d\sigma. 
$$

\par We note that if $g\in {\mathcal V}$,
 $$
 (g,g)=\int_\Stwo(|\nabla g|^2-2g^2)d\sigma
 $$
 is nonnegative and vanishes if and only if $g\in \Pee$. We will quantify this in the lemma below. Thus, it is natural to consider the quotient space
$$
{\mathcal V}/\Pee
$$
with quotient inner product defined as 
\be\label{QuotientInnerProduct}
(\bg,\bh):=(g,h).
\ee
Here we are extending our notation $\bg=g+\Pee$ to $g\in {\mathcal V}$.

\par We claim that ${\mathcal V}/\Pee$ is a Hilbert space. First, we'll derive a basic stability estimate. 
\begin{lem}
For each $g\in {\mathcal V}$, the following inequality holds
\be\label{SpaceVeeStable}
4\min_{a\in \R^3}\int_{\Stwo}(g(u)-a\cdot u)^2d\sigma(u)\le \int_\Stwo(|\nabla g|^2-2g^2)d\sigma.
\ee
\end{lem}
\begin{proof}
Recall the collection of spherical harmonics $\phi^m_\ell$, where $\ell=0,1,2,\dots$ and $m=-\ell,-\ell+1,\dots, \ell-1,\ell$. These are eigenfunctions for $-\Delta$ which form an orthonormal basis for $L^2(\Stwo)$. That is,
$$
\int_{\Stwo}\phi^m_\ell\phi^{m'}_{\ell'}d\sigma=\delta_{mm'}\delta_{\ell\ell'}
$$
and 
$$
-\Delta\phi^m_\ell=\ell(\ell+1)\phi^m_\ell
$$
in $\Stwo$ for each $\ell,\ell'=0,1,2,\dots$, $m=-\ell,-\ell+1,\dots, \ell-1,\ell$, and $m'=-\ell',-\ell'+1,\dots, \ell'-1,\ell'$.  We especially note that $\phi^0_0$ is a nonzero constant and we may assume
$$
\phi^1_{-1}(u)=cu_1,\quad\phi^1_0(u)=cu_2,\quad \phi^1_1(u)=cu_3.
$$
Here $c>0$ is the normalization constant to ensure these functions have $L^2(\Stwo)$ norm equal to 1. 

\par As a result, we may write  
$$
g=\sum^\infty_{\ell=0}\sum^\ell_{m=-\ell}c^m_\ell\phi^m_\ell
$$
where 
$$
c^m_\ell=\int_{\Stwo}g\phi^m_\ell d\sigma.
$$
We note that $c^0_0=0$ since the average of $g$ is $0$. 
It follows from these observations that 
\begin{align}
\int_{\Stwo}g^2d\sigma=\sum^\infty_{\ell=1}\sum^\ell_{m=-\ell}(c^m_\ell)^2
\end{align}
and 
\begin{align}
\int_{\Stwo}|\nabla g|^2d\sigma=\sum^\infty_{\ell=1}\sum^\ell_{m=-\ell}\ell(\ell+1)(c^m_\ell)^2.
\end{align}

\par This leads to 
\begin{align}
\int_{\Stwo}\left(|\nabla g|^2-2g^2\right)d\sigma&=\sum^\infty_{\ell=1}\sum^\ell_{m=-\ell}[\ell(\ell+1)-2](c^m_\ell)^2\\
&=\sum^\infty_{\ell=2}\sum^\ell_{m=-\ell}[\ell(\ell+1)-2](c^m_\ell)^2\\
&\ge 4\sum^\infty_{\ell=2}\sum^\ell_{m=-\ell}(c^m_\ell)^2\\
&= 4\int_{\Stwo}(g(u)-c(c^1_{-1}u_1+c^1_{0}u_2+c^1_{1}u_3))^2d\sigma(u)\\
&= 4\min_{a\in \R^3}\int_{\Stwo}(g(u)-a\cdot u)^2d\sigma(u).
\end{align}
\end{proof}
\begin{cor}
${\mathcal V}/\Pee$ endowed with the inner product $(\cdot,\cdot)$ is complete.
\end{cor}
\begin{proof}
Suppose $(\bg^k)_{k\in \N}\subset {\mathcal V}/\Pee$ is a Cauchy sequence. For each $g^k\in {\mathcal V}$, we select $a^k\in \R^3$ to satisfy
$$
\min_{a\in \R^3}\int_{\Stwo}(g^k(u)-a\cdot u)^2d\sigma(u)=\int_{\Stwo}(g^k(u)-a^k\cdot u)^2d\sigma(u).
$$
We also set $\tilde g^k(u)=g^k(u)-a^k\cdot u$ and note $\tilde g^k\in \bg^k$ for each $k\in \N$. In view of the previous lemma, for each $\epsilon>0$, there is $N\in N$ so that 
$$
\epsilon\ge \int_\Stwo(|\nabla (\tilde g^k-\tilde g^\ell)|^2-2(\tilde g^k-\tilde g^\ell)^2)d\sigma\ge 4\int_{\Stwo}(\tilde g^k-\tilde g^\ell)^2d\sigma
$$
for each $k,\ell\ge N$. As a result, $(\tilde g^k)_{k\in \N}$ is Cauchy in $L^2(\Stwo)$. It also follows from the above inequality that $(\tilde g^k)_{k\in \N}$ is 
 Cauchy in $H^1(\Stwo)$.  Since $H^1(\Stwo)$ is complete and ${\mathcal V}$ is closed, there is $g\in{\mathcal V}$ such that  $\tilde g^k\rightarrow g$ in $H^1(\Stwo)$. We conclude 
$$
(\bg^k-\bg,\bg^k-\bg)=\int_\Stwo(|\nabla (\tilde g^k-g)|^2-2(\tilde g^k-g)^2)d\sigma\rightarrow 0
$$
as $k\rightarrow\infty$.
\end{proof}

\par Suppose $\bg: [0,\infty)\rightarrow {\mathcal V}/\Pee$ is a measurable mapping.  We may consider the integral 
$$
\int^\infty_0(\bg(t),\bg(t))dt=\int^\infty_0\int_{\Stwo}\left(|\nabla g(t)|^2-2g(t)^2\right)d\sigma dt
$$
which may or may not be finite. Whenever this integral is finite, we say that 
$$
\bg \in L^2([0,\infty);{\mathcal V}/\Pee).
$$
Further, we may view this space as a Hilbert space with the inner product 
$$
(\bg_1,\bg_2)\mapsto \int^\infty_0(\bg_1(t),\bg_2(t))dt
$$
if we identify paths that are equal almost everywhere. 

\par A basic assertion that we will make use of is as follows. 
\begin{lem}\label{WeakConvCeeLem}
The collection of $\bg \in L^2([0,\infty);{\mathcal V}/\Pee); t\mapsto \bg(t)$ such that
$$
g(t)\in \Cee
$$
for almost every $t\ge 0$ is weakly closed. 
\end{lem}
\begin{proof}
Let us call the collection of maps in question ${\mathcal K}$. We claim that ${\mathcal K}$ is convex and closed. The result would then follow from Mazur's theorem (Chapter V of \cite{MR617913}). The convexity of ${\mathcal K}$ follows easily since $\Cee$ is convex.  Now suppose $(\bg^k)_{k\in \N}\subset {\mathcal K}$ converges in $L^2([0,\infty); {\mathcal V}/\Pee)$ to some $\bg$. We need to verify $\bg\in {\mathcal K}$. 
\par To this end, we subtract a subsequence $(\bg^{k_j})_{j\in \N}$ such that $\bg^{k_j}(t)\rightarrow \bg(t)$ for almost every $t\ge 0$ in ${\mathcal V}/\Pee$. For any such $t$, Lemma \ref{CompactCeeLem} implies the existence of a sequence $(g^{j}(t))_{j\in \N}\subset \Cee$ with $g^{j}(t)\in \bg^{k_j}(t)$ that converges in $C^1(\Stwo)$ to some limit function $\tilde g(t)\in\Cee$.  Then $( g^{j}(t),h)$ converges to both $( g(t),h)$ and $(\tilde g(t),h)$ 
 for $\bh\in {\mathcal V}/\Pee$. That is, 
 $$
 ( g(t),h)=(\tilde g(t),h)
 $$
 for all $\bh\in {\mathcal V}/\Pee$.  It follows that $g(t)-\tilde g(t)\in \Pee$ which implies $g(t)\in\Cee$. We conclude that for almost every $t\ge 0$, $g(t)\in \Cee$. Thus, ${\mathcal K}$ is closed. 
\end{proof}
The following claim will also be useful for us.
\begin{lem}\label{AchOneMeaswithCeegivesMeas}
Suppose $\bg \in L^2([0,\infty);{\mathcal V}/\Pee); t\mapsto \bg(t)$ satisfies
$$
g(t)\in \Cee
$$
for almost every $t\ge 0$. There is a measurable $$\tilde\bg: [0,\infty)\rightarrow C(\Stwo)/\Pee;t\mapsto \tilde\bg(t)$$ with  $\tilde g(t)\in\Cee$ for $t\ge 0$ and
$$
\tilde\bg(t)=\bg(t)
$$
for almost every $t\ge 0$.
\end{lem}
\begin{proof}
Let $N\subset [0,\infty)$ be a null set for which $g(t)\in \Cee$ for $t\not\in N$. Set 
$$
\tilde\bg(t):=
\begin{cases}
\bg(t), \quad &t\not\in N\\
0,\quad &t\in N.
\end{cases}
$$
Then $\tilde g(t)\in \Cee$ for each $t\ge 0$ and $\tilde\bg(t)=\bg(t)$ for almost every $t\ge 0$. 

\par  Now assume $\mu\in\Pee^\perp$. According to Proposition \ref{smoothlemma} in the appendix: for each $\epsilon\in (0,1)$, there is 
$\mu^\epsilon\in C^\infty(\Stwo)$ 
 such that 
\be\label{myooEspOrtho}
 \int_{\Stwo}u_i \mu^\epsilon d\sigma=0
\ee
 for $i=1,2,3$, and 
$$
\left|\int_\Stwo g(t)\mu^\epsilon d\sigma-\int_\Stwo g(t)d\mu\right|\le  4\sqrt{2\epsilon}\;\|\mu\|_*
$$
for all $t\not\in N$. Therefore, 
\be\label{mugeeteelimit}
\int_\Stwo g(t)d\mu=\lim_{\epsilon\rightarrow 0^+}\int_\Stwo g(t)\mu^\epsilon d\sigma
\ee
for $t\not\in N$.

\par As a result, it suffices to show 
\be\label{myooEspIntFun}
[0,\infty)\ni t\mapsto \int_\Stwo g(t)\mu^\epsilon d\sigma
\ee
is measurable as
$$
[0,\infty)\ni t\mapsto \int_\Stwo g(t)d\mu
$$
would be the almost everywhere limit of measurable functions. Pettis' theorem (Chapter V section 4 of \cite{MR617913}) would then imply that $\tilde \bg: [0,\infty)\rightarrow C(\Stwo)/\Pee$ is measurable. 

\par We now focus on showing \eqref{myooEspIntFun} is measurable.  To this end, we fix an element $\bh\in  {\mathcal V}/\Pee$.  
 In view of \eqref{myooEspOrtho}, we have 
 \begin{align}
 \int_{\Stwo}\mu^\epsilon(u)h(u)d\sigma(u)&= \int_{\Stwo}\mu^\epsilon(u)( h(u)-a\cdot u)d\sigma(u)\\
 &\le\left(\int_{\Stwo}(\mu^\epsilon)^2d\sigma\right)^{1/2}\left(\int_\Stwo(h(u)-a\cdot u)^2d\sigma(u)\right)^{1/2}
 \end{align}
for any $a\in \R^3$. Employing \eqref{SpaceVeeStable}, we additionally find 
 \begin{align}
 \int_{\Stwo}\mu^\epsilon h d\sigma &\le\left(\int_{\Stwo}(\mu^\epsilon)^2d\sigma\right)^{1/2}\frac{1}{2}\left(\int_\Stwo\left(|\nabla h|^2-2h^2\right)d\sigma\right)^{1/2}\\
 &=\frac{1}{2}\left(\int_{\Stwo}(\mu^\epsilon)^2d\sigma\right)^{1/2}(\bh,\bh)^{1/2}.
 \end{align}
 Therefore, the linear functional
 $$
{\mathcal V}/\Pee\ni \bh\mapsto  \int_{\Stwo}\mu^\epsilon h d\sigma
 $$
 is continuous. Since $\bg: [0,\infty)\rightarrow {\mathcal V}/\Pee$ is measureable, it follows that \eqref{myooEspIntFun} is measurable.
\end{proof}
\subsection{The metric induced by $\chi_\Cee^*$}\label{subsecCeePerp}
Let us define
$$
\Cee^\perp:=\{\xi\in M(\Stwo): \xi|_{\Cee}=0 \},
$$
and observe that  $\Cee^\perp\subset \Pee^\perp$ is a closed subspace. We will also consider the quotient space 
$$
\Pee^\perp/{\Cee}^\perp=\{\xi+\Cee^\perp:  \xi\in \Pee^\perp\}.
$$
As we did to ease notation when expressing elements of the quotient space $C(\Stwo)/\Pee$, we will write
$$
\bxi =\xi+\Cee^\perp
$$
whenever $\xi\in \bxi$. 

\par While $\Pee^\perp/{\Cee}^\perp$ admits the standard quotient norm 
$$
\|\bxi\|_{*}:=\inf_{\eta \in\Cee^\perp }\|\xi+\eta\|_{*},
$$
is it also is endowed with another norm 
\be\label{RightDualNorm}
\chi_{\Cee}^*(\bxi):=\sup\left\{\langle \xi,g\rangle: g\in \Cee\right\}.
\ee
We have labeled this norm $\chi_{\Cee}^*$ as it arises as the convex dual of the characteristic function 
$$
\chi_{\Cee}(g)
:=\begin{cases}
0,\quad & g\in \Cee\\
+\infty, \quad & g\not\in\Cee.
\end{cases}
$$
 Indeed the right hand side of \eqref{RightDualNorm} is the dual of $\chi_{\Cee}$; since 
this function is invariant under translations by elements of $\Cee^\perp$, we naturally consider it as a function on $\Pee^\perp/\Cee^\perp$. 

\par First we note that the quotient norm on $\Pee^\perp/\Cee^\perp$ controls the norm $\chi_{\Cee}^*$. In particular, the topology associated with 
$\chi_{\Cee}^*$ is weaker than the topology determined by the quotient norm.  
\begin{prop}
For $\bxi\in \Pee^\perp/\Cee^\perp$, 
\be\label{BestConstQuestEst}
 \chi_{\Cee}^*(\bxi)\le \left(\sqrt{\frac{3}{8}}-\frac{1}{2}\right)\|\bxi\|_{*}.
\ee
\end{prop}
\begin{proof}
We recall that if $g\in \Cee$, then $\|\bg\|\le \sqrt{\frac{3}{8}}-\frac{1}{2}$.  Therefore, 
\begin{align}
\chi_{\Cee}^*(\bxi)&=\sup\left\{\langle \xi,g\rangle: g\in \Cee\right\}\\
&\le \sup\left\{\langle \xi,g\rangle: \|\bg\|\le\sqrt{\frac{3}{8}}-\frac{1}{2}\right\}\\
&=\left( \sqrt{\frac{3}{8}}-\frac{1}{2}\right)\|\bxi\|_{*}.
\end{align}
\end{proof}
Next, we will argue that quotient norm bounded subsets of $\Pee^\perp/\Cee^\perp$ are compact metric spaces when endowed with the metric 
\be
(\bxi,\bzeta)\mapsto \chi_{\Cee}^*(\bxi-\bzeta)
\ee
induced by $\chi_{\Cee}^*$.
\begin{prop}\label{BeeArrCompactPoverC}
For each $r>0$,
$$
{\mathcal B}_r:=\{\bxi\in \Pee^\perp/\Cee^\perp: \|\bxi\|_{*}\le r\}
$$
is a compact metric space when endowed with metric induced by $\chi_{\Cee}^*$.
\end{prop}
\begin{proof}
Suppose that $(\bxi^k)_{k\in \N}\subset {\mathcal B}_r$.  For each $k\in \N$, we may choose $\eta^k\in \Cee^\perp$ such that 
$$
r\ge \|\bxi^k\|_{*}\ge \|\xi^k+\eta^k\|_{*}-\frac{r}{2}.
$$
Let us set $\zeta^k:=\xi^k+\eta^k$ so that 
$$
\|\zeta^k\|_*\le 3r/2 .
$$
By Alaoglu's theorem, that there is a subsequence $(\zeta^{k_j})_{j\in \N}$ that converges weak* in $M(\Stwo)$ to some $\xi\in \Pee^\perp$. 

\par A routine application of Lemma \ref{CompactCeeLem} implies that for each $j\in \N$, there is $g^j\in \Cee$ with 
$$
\chi_{\Cee}^*(\bxi^{k_j}-\bxi)=\langle \xi^{k_j}-\xi,g^j\rangle=\langle \zeta^{k_j}-\xi,g^j\rangle.
$$
We may also assume without any loss of generality that $g^j$ converges uniformly to some $g\in \Cee$. Notice that
\begin{align}
\chi_{\Cee}^*(\bxi^{k_j}-\bxi)&=\langle \zeta^{k_j}-\xi,g^j\rangle\\ 
&=\langle \zeta^{k_j}-\xi,g^j-g\rangle+\langle \zeta^{k_j}-\xi,g\rangle\\ 
&\le \|\zeta^{k_j}-\xi\|_{*} \|g^j-g\| +\langle \zeta^{k_j}-\xi,g\rangle\\
&\le 3r\|g^j-g\|+\langle \zeta^{k_j}-\xi,g\rangle.
\end{align}
Therefore, 
$$
\limsup_{j\rightarrow\infty}\chi_{\Cee}^*(\bxi^{k_j}-\bxi)=0.
$$
\end{proof}
\begin{cor}\label{QuotientNormPeePerpLSC}
Suppose $(\bxi^k)_{k\in \N}\subset \Pee^\perp/\Cee^\perp$ and $\bxi\in\Pee^\perp/\Cee^\perp$ with
\be\label{xikgoestoxiAssump}
\lim_{k\rightarrow\infty}\chi_\Cee^*(\bxi^k-\bxi)=0.
\ee
Then 
$$
 \|\bxi\|_{*}\le \liminf_{k\rightarrow \infty} \|\bxi^k\|_{*}.
$$
\end{cor}
\begin{proof}
Without any loss of generality, suppose 
$$
 \liminf_{k\rightarrow \infty} \|\bxi^k\|_{*}= \lim_{j\rightarrow \infty} \|\bxi^{k_j}\|_{*}:=L<\infty.
$$
For a given $\epsilon>0$, we can choose $N$ such that 
$$
\|\bxi^{k_j}\|_{*}\le L+\epsilon
$$
for $j\ge N$. By the previous proposition, 
$$
\{\bzeta\in\Pee^\perp/\Cee^\perp : \|\bzeta\|_{*}\le L+\epsilon\}
$$
is a compact metric space when endowed with the metric induced by $\chi_\Cee^*$. By assumption \eqref{xikgoestoxiAssump}, it must 
be that 
$$
 \|\bxi\|_{*}\le L+\epsilon.
$$
We conclude since $\epsilon>0$ is arbitrary. 
\end{proof}

\subsection{$E^*$ variation}
Recall that the convex conjugate of $E$ is given by
$$
E^*(\zeta)=\sup\{\langle \zeta,g\rangle -E(g):g\in \Cee \}
$$
for $\zeta\in \Pee^\perp. $  Note that $E^*$ is convex, proper, and weak* lower-semicontinuous. Moreover,  For a given path $\xi: [0,\infty)\rightarrow \Pee^\perp$, we will consider the integrals
$$
\int^t_sE^*(\dot\xi(\tau))d\tau
$$
for $0\le s\le t<\infty$. However, we do not always want to require the almost everywhere weak* differentiability of $\xi$.  This leads us to the following definition.

\begin{defn}
Let $\xi: [0,\infty)\rightarrow \Pee^\perp$ and $0\le s< t<\infty$. The {\it $E^*$ variation} of $\xi$ on $[s,t]$ is defined as 
\be
E^*V(\xi,s,t):=\sup\left\{\sum^N_{k=1}(\tau_k-\tau_{k-1})E^*\left(\frac{\xi(\tau_k)-\xi(\tau_{k-1})}{\tau_k-\tau_{k-1}}\right): s=\tau_0<\dots<\tau_N=t\right\}.
\ee
When $E^*V(\xi,s,t)<\infty$, $\xi$ has {\it finite $E^*$ variation} on $[s,t]$; if this is the case for all $[s,t]\subset [0,\infty)$, $\xi$ has {\it locally finite $E^*$ variation}.  The {\it $E^*$ variation} of $\xi$ on $[0,\infty)$ is defined as 
 $$
 E^*V(\xi,0,\infty):=\lim_{t\rightarrow\infty}E^*V(\xi,0,t)
 $$
provided this limit is finite; in this case, we say $\xi$ has {\it finite $E^*$ variation}.
\end{defn}
An elementary identity is as follows.
\begin{lem}\label{EstarVadditive}
Suppose  $\xi: [0,\infty)\rightarrow \Pee^\perp$ has locally finite $E^*$ variation. Then 
\be\label{MontoneEstarV}
E^*V(\xi,s,t)= E^*V(\xi,0,t)-E^*V(\xi,0,s)
\ee
for all $0\le s< t<\infty$.
\end{lem}
\begin{proof}
Let $\epsilon>0$, choose a partition $0=\tau_0<\dots<\tau_N=s$ such that 
$$
E^*V(\xi,0,s)\le \sum^N_{k=1}(\tau_k-\tau_{k-1})E^*\left(\frac{\xi(\tau_k)-\xi(\tau_{k-1})}{\tau_k-\tau_{k-1}}\right)+\frac{\epsilon}{2},
$$
and select another partition $s=\tau_N<\dots<\tau_M=t$ such that 
$$
E^*V(\xi,s,t)\le \sum^M_{k=N+1}(\tau_k-\tau_{k-1})E^*\left(\frac{\xi(\tau_k)-\xi(\tau_{k-1})}{\tau_k-\tau_{k-1}}\right)+\frac{\epsilon}{2}.
$$
As $0=\tau_0<\dots<\tau_M=t$,  
$$
E^*V(\xi,0,s)+E^*V(\xi,s,t)-\epsilon\le\sum^M_{k=1}(\tau_k-\tau_{k-1})E^*\left(\frac{\xi(\tau_k)-\xi(\tau_{k-1})}{\tau_k-\tau_{k-1}}\right)\le E^*V(\xi,0,t).
$$
Therefore, $E^*V(\xi,0,s)+E^*V(\xi,s,t)\le E^*V(\xi,0,t)$.

\par Again let $\epsilon>0$ and select any partition $0=\tau_0<\dots<\tau_N=t$ of $[0,t]$ with
$$
E^*V(\xi,0,t)\le \sum^N_{j=1}(\tau_j-\tau_{j-1})E^*\left(\frac{\xi(\tau_j)-\xi(\tau_{j-1})}{\tau_j-\tau_{j-1}}\right)+\epsilon.
$$
Suppose there is some $1\le k\le N$ for which $\tau_{k-1}< s<\tau_k$.  Since 
$$
\frac{\xi(\tau_k)-\xi(\tau_{k-1})}{\tau_k-\tau_{k-1}}=\frac{s-\tau_{k-1}}{\tau_k-\tau_{k-1}}\frac{\xi(s)-\xi(\tau_{k-1})}{s-\tau_{k-1}}+\frac{\tau_{k}-s}{\tau_k-\tau_{k-1}}\frac{\xi(\tau_k)-\xi(s)}{\tau_k-s}
$$
and $E^*$ is convex, 
$$
E^*\left(\frac{\xi(\tau_k)-\xi(\tau_{k-1})}{\tau_k-\tau_{k-1}}\right)\le 
\frac{s-\tau_{k-1}}{\tau_k-\tau_{k-1}}E^*\left(\frac{\xi(s)-\xi(\tau_{k-1})}{s-\tau_{k-1}}\right)+\frac{\tau_{k}-s}{\tau_k-\tau_{k-1}}E^*\left(\frac{\xi(\tau_k)-\xi(s)}{\tau_k-s}\right).
$$
It follows that 
\begin{align*}
E^*V(\xi,0,t)&\le \sum^{k-1}_{j=1}(\tau_j-\tau_{j-1})E^*\left(\frac{\xi(\tau_j)-\xi(\tau_{j-1})}{\tau_j-\tau_{j-1}}\right)+
(s-\tau_{k-1})E^*\left(\frac{\xi(s)-\xi(\tau_{k-1})}{s-\tau_{k-1}}\right)\\
&+ (\tau_{k}-s)E^*\left(\frac{\xi(\tau_k)-\xi(s)}{\tau_k-s}\right)+\sum^{N}_{j=k+1}(\tau_j-\tau_{j-1})E^*\left(\frac{\xi(\tau_j)-\xi(\tau_{j-1})}{\tau_j-\tau_{j-1}}\right)+\epsilon\\
&\le E^*V(\xi,0,s)+E^*V(\xi,s,t)+\epsilon,
\end{align*}
where the sum from $j=1,\dots, k-1$ is only present when $k\ge 2$. If $s=\tau_{k-1}$ for some $1\le k\le N$, we can arrive at the same inequality without introducing $s$ as another point in the partition. Therefore, $E^*V(\xi,0,t)\le E^*V(\xi,0,s)+E^*V(\xi,s,t)+\epsilon$ in all cases.  
\end{proof}
We now recall that since $\Pee^\perp$ is the dual of $C(\Stwo)/\Pee$, $\xi:[0,\infty)\rightarrow\Pee^\perp$ is absolutely continuous if and only if $\xi$ weak* differentiable almost everywhere, $\dot\xi\in L^1_{\text{loc}}([0,\infty);\Pee^\perp)$, and 
\be\label{FTC}
\xi(t)-\xi(s)=\int^t_s\dot\xi(\tau)d\tau
\ee
for all $0\le s< t$ (Remark 1.1.3 of \cite{MR2401600}). We can use this characterization to show that if $\xi$ is absolutely continuous then $E^*V(\xi,s,t)$ coincides with the integral that motivated it. 
\begin{prop}\label{classicalVariationNewVariation}
Suppose $\xi: [0,\infty)\rightarrow \Pee^\perp$ is absolutely continuous. Then 
$$
\int^t_sE^*(\dot\xi(\tau))d\tau=E^*V(\xi,s,t)
$$
for all $0\le s\le t$.
\end{prop}
\begin{proof}
Let $0\le s< t$, and suppose $ s=\tau_0<\dots<\tau_N=t$ is a partition of $[s,t]$. In view of \eqref{FTC} and Jensen's inequality 
\begin{align}
(\tau_k-\tau_{k-1})E^*\left(\frac{\xi(\tau_k)-\xi(\tau_{k-1})}{\tau_k-\tau_{k-1}}\right)&=
(\tau_k-\tau_{k-1})E^*\left(\frac{\int^{\tau_k}_{\tau_{k-1}}\dot\xi(\tau)d\tau}{\tau_k-\tau_{k-1}}\right)\\
&\le \int^{\tau_k}_{\tau_{k-1}}E^*(\dot\xi(\tau))d\tau
\end{align}
for $k=1,\dots, N$. Therefore,  
$$
\sum^N_{k=1}(\tau_k-\tau_{k-1})E^*\left(\frac{\xi(\tau_k)-\xi(\tau_{k-1})}{\tau_k-\tau_{k-1}}\right)\le \int^t_sE^*(\dot\xi(\tau))d\tau.
$$
It follows that
$$
E^*V(\xi,s,t)\le \int^{t}_{s}E^*(\dot\xi(\tau))d\tau.
$$

\par Now set 
$$
F(t):=E^*V(\xi,0,t)
$$
for $t\ge 0$. In view of the definition of $E^*$ variation and \eqref{MontoneEstarV}, 
$$
\delta E^*\left(\frac{\xi(t+\delta)-\xi(t)}{\delta}\right)\le E^*V(\xi,t,t+\delta)= F(t+\delta)-F(t)
$$
for $t\ge 0$ and $\delta>0$. As $F$ is monotone and $\xi$ is weak* differentiable almost everywhere, we have 
$$
E^*(\dot\xi(\tau))\le \dot F(t)
$$
for almost every $t\ge 0$.  Here we used that $E^*$ is weak* lower-semicontinuous. That is, 
$$
 \int^{t}_{s}E^*(\dot\xi(\tau))d\tau\le  \int^{t}_{s}\dot F(t)d\tau\le F(t)-F(s)=E^*V(\xi,s,t)
 $$
 for $0\le s\le t$.
\end{proof}
\begin{cor}\label{PiecewiseVariation}
Suppose $\tau>0$ and $(\zeta^k)_{k\ge 0}\subset \Pee^\perp$.  For $t\ge 0$,  define 
$$
\xi(t):=\zeta^{k-1}+\frac{t-(k-1)\tau}{\tau}(\zeta^k-\zeta^{k-1}), \quad (\tau-1)k\le t\le \tau k.
$$
Then
\be
E^*V(\xi,0,N\tau)=\sum^N_{k=1}\tau E^*\left(\frac{\zeta^k-\zeta^{k-1}}{\tau}\right)
\ee
for each $N\in \N$.
\end{cor}
\begin{proof}
Observe that  
$$
\dot\xi(t)=\frac{\zeta^k-\zeta^{k-1}}{\tau}
$$
for $t\in (\tau(k-1), \tau k)$, and  
$$
\int^{\tau k}_{\tau(k-1)}\|\dot\xi(t)\|_{*}dt=\|\zeta^k-\zeta^{k-1}\|_{*}<\infty
$$  
for each $k\in \N$. Furthermore, since $\xi$ is piecewise linear, the fundamental theorem of calculus \eqref{FTC} holds. As a result, $\xi$ is absolutely continuous.  By the previous proposition, 
$$
E^*V(\xi,0,N\tau)= \int^{N\tau}_{0}E^*(\dot\xi(\tau))d\tau=\sum^N_{k=1}\int^{\tau k}_{\tau(k-1)}E^*(\dot\xi(\tau))d\tau=\sum^N_{k=1}\tau E^*\left(\frac{\zeta^k-\zeta^{k-1}}{\tau}\right).
$$\end{proof}
As $E^*(\zeta_1)=E^*(\zeta_2)$ whenever $\zeta_1-\zeta_2\in \Cee^\perp$, $E^*$ can be viewed as a function on the quotient space 
$\Pee^\perp/\Cee^\perp$:
$$
E^*(\bzeta):=\sup\{\langle \zeta,g\rangle -E(\bg): g\in \Cee \}
$$ 
for $\bzeta\in \Pee^\perp/\Cee^\perp.$ In particular, we may consider the $E^*$ variation of paths $\bxi:[0,\infty)\rightarrow\Pee^\perp/\Cee^\perp$ and each of our results above applies to such paths.  We finally note that this extension of $E^*$ is lower-semicontinuous on $(\Pee^\perp/\Cee^\perp,\chi_\Cee^*)$.

\begin{lem}\label{LSCextensionEstar}
Suppose $(\bxi^k)_{k\in \N}$ is a sequence in $\Pee^\perp/\Cee^\perp$ with 
$$
\lim_{k\rightarrow\infty}\chi_\Cee^*(\bxi^k-\bxi)=0
$$
Then 
$$
\liminf_{k\rightarrow\infty}E^*(\bxi^k)\ge E^*(\bxi).
$$
\end{lem}
\begin{proof}
Let $g\in \Cee$. Observe that 
$$
|\langle \xi^k,g\rangle-\langle \xi,g\rangle|=|\langle \xi^k-\xi,g\rangle|\le \chi_\Cee^*(\bxi^k-\bxi).
$$
Therefore, $\langle \xi^k,g\rangle\rightarrow \langle \xi,g\rangle$ as $k\rightarrow\infty$. This implies 
$$
\liminf_{k\rightarrow\infty}E^*(\bxi^k)\ge \liminf_{k\rightarrow\infty} (\langle \xi^k,g\rangle-E(\bg))= \langle \xi,g\rangle-E(\bg).
$$
We conclude upon taking the supremum of $g\in \Cee$.
\end{proof}

\section{Existence theorem}\label{ExistenceSect}
We will now analyze solutions $\xi:[0,\infty)\rightarrow \Pee^\perp$ of the equation 
\be\label{DNExi2}
\partial E^*(\dot\xi(t))+{\mathcal J}^*(\xi(t))\ni 0\quad \text{a.e. $t\ge 0$}.
\ee
These types of flows have been considered in models for diverse phenomena such phase transitions \cite{MR1423808,MR983686,MR1769184}, fracture mechanics \cite{MR2244798,MR2525111,MR2434063}, and hysteresis effects \cite{MR2290410,MR1329094}.  The existence of solutions to doubly nonlinear evolutions 
in Hilbert spaces \cite{MR380532,MR539774, MR1070845} and reflexive Banach spaces \cite{MR1170721} for given initial conditions were established many years ago. 
The reflexivity requirement can be relaxed using the notion of ``curves of maximal slope" as described in Chapter 1 of \cite{MR2401600}.  

\par All of these results require that the corresponding $E^*$ is coercive or superlinear with respect to the norm; this would ensure that solutions are absolutely continuous. 
 We cannot use any of these results directly as our $E^*$ functional is not coercive. A similar problem has been encountered in the theory of rate-independent doubly nonlinear evolutions  \cite{MR4158534,MR4197283,MR2210284,MR2105969,MR2290410,MR2887927,MR3531671,MR3740380,MR3636535}.  Here the corresponding $E^*$ function is a norm and solutions are typically considered as mappings of bounded variation measured in this norm. We will adapt this approach to our problem using the $E^*$ variation notion discussed in the previous section.

\par We will define a type of weak solution of \eqref{DNExi2} below. To avoid confusion, we'll say that $\xi$ is a {\it classical solution} of  \eqref{DNExi2} if it satisfies the equation as described in the introduction. That is, $\xi$ is a classical solution if it is absolutely continuous and if there is a measurable $\bg: [0,\infty)\rightarrow C(\Stwo)/\Pee$ such that
\be\label{geeacheconditions2}
\bg(t)\in {\mathcal J}^*(\xi(t))\cap (-\partial E^*(\dot\xi(t)))
\ee
for almost every $t\ge 0$. Here $\dot\xi(t)$ is the weak* derivative of $\xi$ at a time $t\ge 0$.

\begin{defn}\label{WeakSolutionDefn}
A measurable mapping
$$
\xi: [0,\infty)\rightarrow\Pee^\perp
$$
with locally finite $E^*$ variation is a {\it weak solution} of the doubly nonlinear evolution \eqref{DNExi2} provided there is
a measurable
$$
\bg: [0,\infty)\rightarrow C(\Stwo)/\Pee
$$
such that 
\be\label{DNEChainRuleDef}
\frac{1}{2}\|\xi(s)\|_*^2- \frac{1}{2}\|\xi(t)\|_*^2= E^*V(\xi,s,t)+\int^t_sE(\bg(\tau))d\tau
\ee
for almost every $0\le s\le t$ and 
$$
\bg(t)\in \partial {\mathcal J^*}(\xi(t))
$$
for almost every $t\ge 0$.
\end{defn}

\par For any weak solution $\xi$, $\|\xi(t)\|_*$ is a nonincreasing function of $t$ outside of a null set. We will say such functions are {\it essentially nonincreasing}. In the appendix, we will recall that an essentially nonincreasing function is simply a nonincreasing function modified on a null set (Lemma \ref{EssNonIncLemma}). Also note that by the duality formula \eqref{DualityFormula}
$$
[0,\infty)\ni t\mapsto \|\bg(t)\|
$$
is essentially nonincreasing, as well. 

\par Let us check that the notion of weak solutions extends the classical notion of solution.
\begin{lem}
Assume $\xi: [0,\infty)\rightarrow\Pee^\perp$ is a classical solution of \eqref{DNExi2}. Then $\xi$ is a weak solution of \eqref{DNExi2}.
\end{lem}
\begin{proof}
Since $\xi$ is absolutely continuous, it is measurable. We may also
select a measurable $\bg: [0,\infty)\rightarrow C(\Stwo)/\Pee$  such that \eqref{geeacheconditions2} holds almost everywhere on $[0,\infty)$. 
Note that for almost every $t>0$, 
$$
\frac{1}{2}\|\xi(t+\tau)\|^2_*\ge \frac{1}{2}\|\xi(t)\|^2_*+\langle \xi(t+\tau)-\xi(t),g(t)\rangle 
$$
for $\tau\in \R$ sufficiently small. It follows that 
$$
\lim_{\tau\rightarrow 0^+}\frac{\|\xi(t+\tau)\|^2_*-\|\xi(t)\|^2_*}{2\tau}\ge \langle \dot\xi(t),g(t)\rangle\ge \lim_{\tau\rightarrow 0^-}\frac{\|\xi(t+\tau)\|^2_*-\|\xi(t)\|^2_*}{2\tau}
$$
for almost every $t>0$. As $[0,\infty)\ni t\rightarrow\|\xi(t)\|_*$ is absolutely continuous, 
$$
\frac{d}{dt}\frac{1}{2}\|\xi(t)\|^2_*=\langle \dot\xi(t),g(t)\rangle
$$
almost everywhere in $[0,\infty)$.

\par Recall that $-\dot\xi(t)\in \partial E(\bg(t))\neq \emptyset$ for almost every $t\ge 0$. As a result, $g(t)\in \Cee$ and
$$
\frac{d}{dt}\frac{1}{2}\|\xi(t)\|^2_*=-\langle \dot\xi(t),-g(t)\rangle=-[E^*(\dot\xi(t))+E(\bg(t))]
$$
for almost every $t\ge 0$. Integrating this identity from $s$ to $t$ with $s\le t$ gives
\begin{align}
\frac{1}{2}\|\xi(s)\|_*^2- \frac{1}{2}\|\xi(t)\|_*^2&= \int^t_sE^*(\dot\xi(\tau))d\tau+\int^t_sE(\bg(\tau))d\tau\\
& = E^*V(\xi,s,t)+\int^t_sE(\bg(\tau))d\tau.
\end{align}
Here we used Proposition \ref{classicalVariationNewVariation}. We conclude that that $\xi$ is indeed a weak solution.
\end{proof}

\par The main theorem of this paper as follows. 
\begin{thm}\label{mainThm}
Assume $g^0\in \Cee$ and that $\xi^0\in \Jay(\bg^0)$. There is a weak solution $\xi$ with finite $E^*$ variation which satisfies the initial condition
\be\label{xizeroCond}
\xi(0)|_{\Cee}=\xi^0|_{\Cee}.
\ee
Moreover, there is a corresponding $\bg:[0,\infty)\rightarrow C(\Stwo)/\Pee$ 
as described in Definition \ref{WeakSolutionDefn} for which $\bg\in L^2([0,\infty);{\mathcal V}/\Pee)$ and 
$$
[0,\infty)\ni t\mapsto E(\bg(t))\quad\text{is essentially nonincreasing}.
$$
\end{thm}
As a corollary to the theorem above, we have the $t E(\bg(t))$ and $\|\bg(t)\|$ both converge to zero along sequences to times $t$ tending to infinity outside of 
a null set. 
\begin{cor}
There is a null set $N\subset [0,\infty)$ such that 
\be\label{LargeTimeLimits}
\lim_{\substack{t\rightarrow\infty\\ t\not\in N}}t E(\bg(t))=0\quad \text{and}\quad \lim_{\substack{t\rightarrow\infty\\ t\not\in N}}\|\bg(t)\|=0.
\ee
\end{cor}
\begin{proof}
As $E\circ \bg:[0,\infty)\rightarrow [0,\infty)$ is essentially nonincreasing, Proposition \eqref{teeffteeGoesZero} of the appendix implies there is a null set $N$ for which 
\be\label{OneHalfLargeTimeLimits}
\lim_{\substack{t\rightarrow\infty\\ t\not\in N}}t E(\bg(t))=0.
\ee
\par Without any loss of generality, we may assume $g(t)\in\Cee$ for all $t\not\in N$ since this occurs for almost every time $t\ge 0$. Let $(t_k)_{k\in \N}\subset [0,\infty)\setminus N$ be a sequence increasing to $\infty$. By Lemma \ref{CompactCeeLem}, $(\bg(t_k))_{k\in \N}$ has a convergent subsequence 
$(\bg(t_{k_j}))_{j\in \N}$. Furthermore, there is a sequence $(g^j)_{j\in \N}\subset C^1(\Stwo)$ with $g^j\in \bg(t_{k_j})$ and which converges in $C^1(\Stwo)$ to some limit function $g^\infty$. By \eqref{OneHalfLargeTimeLimits}, 
$$
E(\bg(t_{k_j}))=\int_{\Stwo}\left(\frac{1}{2}|\nabla g^j|^2-(g^j)^2\right)d\sigma\le \frac{1}{t_{k_j}}
$$
for all large enough $j$. Thus, 
$$
\int_{\Stwo}\left(\frac{1}{2}|\nabla g^\infty|^2-(g^\infty)^2\right)d\sigma=\lim_{j\rightarrow\infty}\left(\int_{\Stwo}\frac{1}{2}|\nabla g^j|^2-(g^j)^2\right)d\sigma=0.
$$
This implies $g^\infty(u)=b\cdot u$ for some $b\in \R^3$. It then follows that 
$$
\lim_{j\rightarrow\infty}\|\bg(t_{k_j})\|\le \lim_{j\rightarrow\infty}\max_{|u|=1}|g^j(u)-b\cdot u|=0.
$$
Since the sequence $(t_k)_{k\in \N}\subset [0,\infty)\setminus N$ was arbitrary, we conclude \eqref{LargeTimeLimits}.
\end{proof}
We may interpret these results geometrically. Choose a null set $N\subset [0,\infty)$ such that for all $t\in N^c$, $g(t)\in \Cee$.  For each $t\in N^c$, $g(t)+1/2$ is the $\Stwo$ restriction of the support function of a constant width body $K_t$.  Then Theorem \ref{mainThm}, the limits \eqref{LargeTimeLimits}, and equation \eqref{propFormulafofR2} imply the following assertion on the family $(K_t)_{t\in N^c}$. 
\begin{cor} 
(i) $R(K_t)=\frac{1}{2}+\|\bg(t)\|$ is essentially nonincreasing, and 
$$
\displaystyle\lim_{\substack{t\rightarrow\infty\\ t\not\in N}}R(K_t)=\frac{1}{2}.
$$
(ii)
$V(K_t)=\frac{\pi}{6}-\frac{1}{2}E(\bg(t))$ is essentially nondecreasing, and
$$
\displaystyle\lim_{\substack{t\rightarrow\infty\\ t\not\in N}}t\left(\frac{\pi}{6}-V(K_t)\right)=0.
$$
(iii) Suppose $a(t)\in \R^3$ is the center of the circumball of $K_t$ for each $t\in N^c$. Then 
 $$
\displaystyle\lim_{\substack{t\rightarrow\infty\\ t\not\in N}}d_{{\mathcal H}}(K_t-a(t),B_{1/2}(0))=0
$$
\end{cor}

\par We will design a weak solution as asserted in Theorem \ref{mainThm} as follows.  First, we will show how to solve a discrete version of \eqref{DNExi2} which depends on a parameter $\tau>0$. Next, we will show how to use these solutions to form a family of approximate solutions of \eqref{DNExi2} indexed by $\tau$. Then we will derive various inequalities satisfied by this family of approximate solutions which are independent of $\tau$. Finally, we will explain how to extract a sequence of $\tau$ tending to $0$ for which the approximate solutions converge to a weak solution.  

\subsection{Implicit time scheme}
For the remainder of this section we will suppose the hypotheses of Theorem \ref{mainThm}. That is, we will assume $g^0\in \Cee$ and that $\xi^0\in \Jay(\bg^0)$. In order to prove  this theorem, we will use the following implicit time scheme: fix $\tau>0$ and find a sequence $(\xi^k)_{k\in \N}\subset \Pee^\perp$  such that 
\be\label{ITSxi}
\partial E^*\left(\frac{\xi^{k}-\xi^{k-1}}{\tau}\right)+\Jay^*(\xi^k)\ni 0
\ee
for $k\in \N$.  Let us verify that there is indeed a solution sequence.


\begin{lem}
There is a solution sequence $(\xi^k)_{k\in \N}$ of \eqref{ITSxi} which satisfies 
\be\label{chseekMin}
\|\xi^k\|_*=\|\bxi^k\|_*
\ee
for each $k\in \N$. 
\end{lem}
\begin{rem}
Here $\bxi^k=\xi^k+\Cee^\perp$ is the equivalence class defined in subsection \ref{subsecCeePerp}.  
\end{rem}
\begin{proof}
We will proceed by induction.  Once $\xi^0,\dots, \xi^{k-1}\in \Pee^\perp$ are determined, we can minimize 
\be\label{ITSfunctional}
\Pee^\perp\ni \xi\mapsto \tau E^*\left(\frac{\xi-\xi^{k-1}}{\tau}\right)+\frac{1}{2}\|\xi\|_*^2
\ee
over $\Pee^\perp$.  Starting with any minimizing sequence, we can employ Alaoglu's theorem and the weak* lower-semicontinuity of both $E^*$ and the norm on $\Pee^\perp$ to conclude the existence of a minimizer $\xi^k\in\Pee^\perp$.  Therefore, there is a solution sequence $(\xi^k)_{k\in \N}$ of \eqref{ITSxi}. 

\par We may select 
$$
\bg^k\in \Jay^*(\xi^k)\cap\left(-\partial E^*\left(\frac{\xi^k-\xi^{k-1}}{\tau}\right)\right)
$$
for $k\in \N$.  Note that since 
\be\label{EbgkSubdiff}
-\frac{\xi^k-\xi^{k-1}}{\tau}\in \partial E(\bg^k)\neq \emptyset,
\ee
$g^k\in \Cee$ for each $k\in \N$, as well.  Therefore, if $\eta\in\Cee^\perp$, then 
$$
\frac{1}{2}\|\xi^k+\eta\|^2_*\ge \frac{1}{2}\|\xi^k\|^2_*+\langle \eta,g^k\rangle=\frac{1}{2}\|\xi^k\|^2_*.
$$
We conclude \eqref{chseekMin}.
\end{proof}

\par As we saw for the classical solutions of \eqref{DNExi2}, the sequences  $(\xi^k)_{k\in \N}\subset \Pee^\perp$ and 
$(\bg^k)_{k\in \N}\subset C(\Stwo)/\Pee$ have two important monotonicity properties. 

\begin{lem}
For each $0\le j<k$,
\be\label{DiscreteMonForm1}
\frac{1}{2}\|\xi^{j}\|_*^2\ge \frac{1}{2}\|\xi^{k}\|_*^2 +\sum^{k}_{\ell=j+1}\tau \left[E^*\left(\frac{\xi^{\ell}-\xi^{\ell-1}}{\tau}\right)+E(\bg^\ell)\right]
\ee
and 
\be\label{DiscreteMonForm2}
E(\bg^{j})\ge E(\bg^{k})+\sum^{k}_{\ell=j+1}\left\langle \frac{\xi^{k}-\xi^{k-1}}{\tau},g^{k}-g^{k-1}\right\rangle.
\ee
\end{lem}
\begin{proof} First note
\begin{align}
\frac{1}{2}\|\xi^{\ell-1}\|_*^2& \ge \frac{1}{2}\|\xi^{\ell}\|_*^2 +\langle \xi^{\ell-1}-\xi^\ell,g^\ell\rangle\\
& = \frac{1}{2}\|\xi^{\ell}\|_*^2+\tau\left\langle \frac{\xi^{\ell}-\xi^{\ell-1}}{\tau},-g^\ell\right\rangle\\
&= \frac{1}{2}\|\xi^{\ell}\|_*^2 +\tau \left[E^*\left(\frac{\xi^{\ell}-\xi^{\ell-1}}{\tau}\right)+E(-\bg^\ell)\right]\\
&= \frac{1}{2}\|\xi^{\ell}\|_*^2 +\tau \left[E^*\left(\frac{\xi^{\ell}-\xi^{\ell-1}}{\tau}\right)+E(\bg^\ell)\right].
\end{align}
Then \eqref{DiscreteMonForm1} follows from summing from $\ell=j+1$ to $k$.

\par Next observe that in view of \eqref{EbgkSubdiff}, 
\begin{align}
E(\bg^{\ell-1})&\ge E(\bg^{\ell})+\left\langle - \frac{\xi^{\ell}-\xi^{\ell-1}}{\tau},g^{\ell-1}-g^{\ell}\right\rangle\\
&= E(\bg^{\ell})+\left\langle \frac{\xi^{\ell}-\xi^{\ell-1}}{\tau},g^{\ell}-g^{\ell-1}\right\rangle.
\end{align}
Inequality \eqref{DiscreteMonForm2} results from summing from $\ell=j+1$ to $k$.
\end{proof}
\begin{rem}
 Since 
$$
\|\xi^{k}\|_* = \|\bg^{k}\|,
$$
it follows that $\left( \|\bg^{k}\|\right)_{k\in \N}$ is nonincreasing.    
\end{rem}

\par In designing a solution, it will help to identify some key variables.   We will denote 
$$
\xi_\tau(t):=\xi^{k-1}+\left(\frac{t-(k-1)\tau}{\tau}\right)\left(\xi^{k}-\xi^{k-1}\right)
$$
for $t\in [(k-1)\tau,k\tau]$ and
$$
\zeta_\tau(t)
:=
\begin{cases}
\xi^0, \quad & t=0\\
\xi^k, \quad & t\in  ((k-1)\tau,k\tau]
\end{cases}
$$ 
for $t\ge 0$.  In addition, we will consider 
$$
\bg_\tau(t)
:=
\begin{cases}
\bg^0, \quad & t=0\\
\bg^k, \quad & t\in ((k-1)\tau,k\tau]
\end{cases}
$$  
for $t\ge 0$.  In terms of these variables, we note 
\be
\partial E^*\left(\dot\xi_\tau(t)\right)+\Jay^*(\zeta_\tau(t))\ni 0
\ee
and 
$$
\bg_\tau(t)\in \Jay^*(\zeta_\tau(t))\cap\left(-\partial E^*\left(\dot\xi_\tau(t)\right)\right)
$$
for $t\neq k\tau$.

\subsection{Various bounds}
We will now derive various bounds on the variables we defined. We will see that the following proposition essentially follows 
from the monotonicity formula \eqref{DiscreteMonForm1} and \eqref{DiscreteMonForm2}.

\begin{prop}
For $\tau>0$ and $t\ge 0$, the following inequalities hold.
\be\label{zetatauUpper}
\|\zeta_\tau(t)\|_*\le \|\xi^0\|_*
\ee
\be\label{xitauUpper}
\|\xi_\tau(t)\|_*\le \|\xi^0\|_*
\ee
\be\label{EstarXiTauUpper}
E^*V(\xi_\tau,0,\infty)\le \frac{1}{2} \|\xi^0\|_*^2
\ee
\be\label{EgTauIntUpper}
\int^\infty_0E(\bg_\tau(t))dt\le \frac{1}{2} \|\xi^0\|_*^2
\ee
\be\label{EgTauUpper}
E(\bg_\tau(t))\le E(\bg^0)
\ee
\end{prop}
\begin{proof}
In view of \eqref{DiscreteMonForm1}, $\|\xi^k\|_*\le\|\xi^{k-1}\|_*\le \|\xi^0\|_*$ for $k\in \N$. It is now immediate that \eqref{zetatauUpper} holds. Likewise
for $t\in [(k-1)\tau,k\tau]$, 
\begin{align*}
\|\xi_\tau(t)\|_*&=\left\|\left(1-\frac{t-(k-1)\tau}{\tau}\right)\xi^{k-1}+\left(\frac{t-(k-1)\tau}{\tau}\right)\xi^k\right\|_*\\
&\le \left(1-\frac{t-(k-1)\tau}{\tau}\right)\|\xi^{k-1}\|_*+\left(\frac{t-(k-1)\tau}{\tau}\right)\|\xi^k\|_*\\
&\le \|\xi^{k-1}\|_*\\
&\le \|\xi^0\|_*.
\end{align*}
We conclude \eqref{xitauUpper}. 

\par In view of Corollary \ref{PiecewiseVariation} and \eqref{DiscreteMonForm1}, 
\be
E^*V(\xi_\tau,0,N\tau)=\sum^N_{k=1}\tau E^*\left(\frac{\zeta^k-\zeta^{k-1}}{\tau}\right)\le \frac{1}{2}\|\xi^0\|_*^2
\ee
for each $N\in \N$. Sending $N\rightarrow\infty$ gives \eqref{EstarXiTauUpper}. Likewise, we have 
\be
\int^{N\tau}_0E(\bg_\tau(t))dt=\sum^N_{k=1}E(\bg^k)\tau \le \frac{1}{2}\|\xi^0\|_*^2
\ee
for every $N\in \N$. And sending $N\rightarrow\infty$ gives \eqref{EgTauIntUpper}.

\par By \eqref{DiscreteMonForm2}, $E(\bg^k)\le E(\bg^0)$ for $k\in \N$. Therefore, \eqref{EgTauUpper} holds. 
\end{proof}

\par It is natural to anticipate that $\xi_\tau$ and $\zeta_\tau$ are close. We can measure their closeness via their respective equivalence classes in $\Pee^\perp/\Cee^\perp$.
\begin{prop}\label{L1xizetaEstTau}
For each $T\ge 0$,
$$
\int^T_0\chi_\Cee^*( \bxi_\tau(t)-\bzeta_\tau(t))dt\le \frac{\tau}{2}\left[\frac{1}{2} \|\xi^0\|_*^2+\frac{\pi}{3}(T+\tau)\right].
$$
\end{prop}
\begin{proof} 
Let $g\in \Cee$ and recall that $E(\bg)\le \pi/3$ by \eqref{BasicBoundonE}.  Note that for $t\in ((k-1)\tau,k\tau]$
\begin{align}
\langle \xi_\tau(t)-\zeta_\tau(t),g\rangle&=\tau\left\langle\frac{\xi^{k}-\xi^{k-1}}{\tau},g\right\rangle \left(\frac{t-k\tau}{\tau}\right)\\
&\le \tau\left[E^*\left(\frac{\xi^{k}-\xi^{k-1}}{\tau}\right) +E(\bg)\right] \left(\frac{k\tau-t}{\tau}\right)\\
&\le \tau\left[E^*\left(\frac{\xi^{k}-\xi^{k-1}}{\tau}\right) +\frac{\pi}{3}\right] \left(\frac{k\tau-t}{\tau}\right).
\end{align}
As a result, 
$$
\chi_\Cee^*( \bxi_\tau(t)-\bzeta_\tau(t))\le \tau\left[E^*\left(\frac{\xi^{k}-\xi^{k-1}}{\tau}\right) +\frac{\pi}{3}\right] \left(\frac{k\tau-t}{\tau}\right).
$$
Integrating over $[(k-1)\tau,k\tau]$ gives 

\begin{align}
\int^{k\tau}_{(k-1)\tau}\chi_\Cee^*( \bxi_\tau(t)-\bzeta_\tau(t))dt&\le \frac{1}{2}\tau^2\left[ E^*\left(\frac{\xi^{k}-\xi^{k-1}}{\tau}\right) +\frac{\pi}{3}\right]
\end{align}
for each $k\in \N$. 

\par Now let $T>0$ and choose $N\in \N$ so  that 
$$
N-1<\frac{T}{\tau}\le N.
$$
In view of \eqref{DiscreteMonForm1},
\begin{align}
\int^T_0\chi_\Cee^*( \bxi_\tau(t)-\bzeta_\tau(t))dt&\le \int^{N\tau}_0\chi_\Cee^*( \bxi_\tau(t)-\bzeta_\tau(t))dt\\
&=\sum^N_{k=1}\int^{k\tau }_{(k-1)\tau}\chi_\Cee^*( \bxi_\tau(t)-\bzeta_\tau(t))dt\\
&\le \sum^N_{k=1}\frac{1}{2}\tau^2\left[ E^*\left(\frac{\xi^{k}-\xi^{k-1}}{\tau}\right) +\frac{\pi}{3}\right]\\
&\le \frac{\tau}{2}\left[\frac{1}{2} \|\xi^0\|_*^2+(N\tau)\frac{\pi}{3}\right]\\
&\le \frac{\tau}{2}\left[ \frac{1}{2} \|\xi^0\|_*^2+(T+\tau)\frac{\pi}{3}\right].
\end{align}
\end{proof}

\par A simple continuity estimate that will prove to be very useful to us is as follows. 
\begin{prop}\label{mainEquiContinuityEstimate}
For $0\le s< t$,
$$
\chi_\Cee^*(\bxi_\tau(t)-\bxi_\tau(s))\le E^*V(\xi_\tau,s,t)+\frac{\pi}{3}(t-s).
$$
\end{prop}
\begin{proof}
Let $g\in \Cee$ and observe 
\begin{align}
\langle \xi_\tau(t)-\xi_\tau(s),g\rangle&=(t-s)\left\langle \frac{\xi_\tau(t)-\xi_\tau(s)}{t-s},g\right\rangle\\
&\le (t-s)E^*V\left(\frac{\xi_\tau(t)-\xi_\tau(s)}{t-s}\right)+(t-s)E(\bg)\\
&\le E^*V(\xi_\tau,s,t)+\frac{\pi}{3}(t-s).
\end{align}
We conclude by taking the supremum over $g\in \Cee$. 
\end{proof}

\subsection{Compactness}
In this subsection, we will establish various assertions involving the convergence of a given quantity along a subsequence of $\tau_j\rightarrow 0^+$.
Since there will be only finitely many of these types of statements, we will not alter the subsequence for each additional limiting assertion which may only guarantee that a limit holds upon passing to a further subsequence.  
 
\par In view of \eqref{DiscreteMonForm1} and \eqref{zetatauUpper},
$$
\|\zeta_{\tau}(t)\|_*\le \|\zeta_{\tau}(s)\|_*\le \|\xi^0\|_*
$$
for $0\le s\le t$.  Likewise \eqref{EstarXiTauUpper} implies
$$
E^*V(\xi_{\tau},0,s)\le E^*V(\xi_{\tau},0,t)\le \frac{1}{2}\|\xi^0\|^2_*
$$
for $0\le s\le t$. By Helly's selection theorem, there is a sequence $(\tau_j)_{j\in \N}$ decreasing to $0$ for which  
\be\label{SimpleHellyConclusion}
\begin{cases}
A(t):=\displaystyle\lim_{j\rightarrow\infty}\|\zeta_{\tau_j}(t)\|_*\\\\
B(t):=\displaystyle\lim_{j\rightarrow\infty}E^*V(\xi_{\tau_j},0,t)
\end{cases}
\ee
for all $t\ge 0$.

\par By \eqref{EgTauIntUpper}, $(\bg_\tau)_{\tau>0}\subset L^2([0,\infty);{\mathcal V}/\Pee)$ is bounded.  Therefore, 
passing to a subsequence if necessary, 
$$
\text{$(\bg_{\tau_j})_{j\in \N}$ converges weakly to some $\bg\in L^2([0,\infty);{\mathcal V}/\Pee)$.}
$$
Recall that $g(t)\in \Cee$ for almost all $t\ge 0$ by Lemma \ref{WeakConvCeeLem}. In view of Lemma \ref{AchOneMeaswithCeegivesMeas}, 
we may also assume 
$$
\begin{cases}
\text{$g(t)\in\Cee$ for all $t\ge 0$ and}\\
\text{$\bg: [0,\infty)\rightarrow C(\Stwo)/\Pee$ is measurable}
\end{cases}
$$
without any loss of generality.  We shall do so going forward.  

\par Our most important convergence assertion is as follows.  
\begin{prop}
There is $\bxi: [0,\infty)\rightarrow\Pee^\perp/\Cee^\perp$ with finite $E^*$ variation for which we can pass to a subsequence if necessary 
to obtain
\be\label{bxiLim}
\lim_{j\rightarrow\infty}\chi^*_\Cee(\bxi_{\tau_j}(t)-\bxi(t))=0
\ee
and
\be\label{IntegralLimzetaJay}
\lim_{j\rightarrow\infty}\int^t_0\chi^*_\Cee(\bzeta_{\tau_j}(s)-\bxi(s))ds=0
\ee
for all $t\ge 0$, and
\be\label{bxiCont}
\lim_{s\rightarrow t}\chi^*_\Cee(\bxi(t)-\bxi(s))=0
\ee
for all but countably many $t\ge 0$.
\end{prop}
\begin{proof}
By \eqref{xitauUpper}, 
$$
\|\bxi_\tau(t)\|_*\le \|\xi_\tau(t)\|_*\le \|\xi^0\|_*;
$$
that is, 
$$
\bxi_\tau(t)\in \{\bzeta\in \Pee^\perp/\Cee^\perp: \|\bzeta\|_*\le \|\xi^0\|_*\}:={\mathcal B}_{\|\xi^0\|_*}
$$
for all $t\ge 0$ and $\tau>0$.  By Proposition \ref{BeeArrCompactPoverC}, $({\mathcal B}_{\|\xi^0\|_*}, \chi_\Cee^*)$
is a compact metric space. 

\par We now appeal to Proposition \ref{mainEquiContinuityEstimate} which asserts
\begin{align}
\chi_\Cee^*(\bxi_\tau(t)-\bxi_\tau(s))&\le E^*V(\xi_\tau,s,t)+\frac{\pi}{3}(t-s)\\
&=E^*V(\xi_\tau,0,t)-E^*V(\xi_\tau,0,s)+\frac{\pi}{3}(t-s)
\end{align}
for $0\le s\le t$. For the last equality, we used Lemma \ref{EstarVadditive}. In view of \eqref{SimpleHellyConclusion},   
\be
\limsup_{j\rightarrow\infty}\chi_\Cee^*(\bxi_{\tau_j}(t)-\bxi_{\tau_j}(s))\le B(t)-B(s)+\frac{\pi}{3}(t-s)
\ee
for all $s\le t$. Since $B:[0,\infty)\rightarrow [0,\infty)$ is nondecreasing, it is continuous for all but countably many $t\ge 0$. 
By an abstract version of the Arzel\`a-Ascoli theorem (Proposition 3.3.1 in \cite{MR2401600}), there is a mapping $\bxi: [0,\infty)\rightarrow {\mathcal B}_{\|\xi^0\|_*}$ 
satisfying \eqref{bxiCont} for all but countably many $t\ge 0$. Moreover, along an appropriate subsequence of $(\bxi_{\tau_j})_{j\in \N}$, \eqref{bxiLim} holds for every $t\ge 0$.  

\par By Lemma \ref{LSCextensionEstar}, $E^*$ is lower-semicontinuous with respect to convergence in the metric induced by $\chi_\Cee^*$. In view of the pointwise convergence \eqref{bxiLim},
$$
E^*V(\bxi,0,t)\le \lim_{j\rightarrow\infty}E^*V(\bxi_{\tau_j},0,t)\le \frac{1}{2}\|\xi^0\|_*^2
$$
for $t\ge 0$. Therefore, $\bxi$ has finite $E^*$ variation. 

\par Recall inequality \eqref{BestConstQuestEst}. This implies 
$$
\chi^*_\Cee(\bxi_{\tau_j}(s)-\bxi(s))\le c\|\bxi_{\tau_j}(s)-\bxi(s)\|_*\le c(\|\bxi_{\tau_j}(s)\|+\|\bxi(s)\|_*)\le 2c\|\xi^0\|_*
$$
for each $s\ge 0$, where $c=\sqrt{\frac{3}{8}}-\frac{1}{2}.$  As a result, we can use dominated convergence to conclude
\be\label{DominateConvergencexitaujay}
\lim_{j\rightarrow\infty}\int^t_0\chi^*_\Cee(\bxi_{\tau_j}(s)-\bxi(s))ds=0
\ee
for each $t\ge 0$. 

\par Let us also recall Proposition \ref{L1xizetaEstTau}, which gives   
$$
\int^t_0\chi^*_\Cee(\bxi_{\tau_j}(s)-\bzeta_{\tau_j}(s))ds\le \frac{\tau_j}{2}\left[\frac{1}{2} \|\xi^0\|_*^2+\frac{\pi}{3}(t+\tau_j)\right].
$$
Employing the triangle inequality, 
$$
\int^t_0\chi^*_\Cee(\bzeta_{\tau_j}(s)-\bxi(s))ds\le \frac{\tau_j}{2}\left[\frac{1}{2} \|\xi^0\|_*^2+\frac{\pi}{3}(t+\tau_j)\right]+\int^t_0\chi^*_\Cee(\bxi_{\tau_j}(s)-\bxi(s))ds.
$$
We can now use  \eqref{DominateConvergencexitaujay} and send $j\rightarrow\infty$ to deduce \eqref{IntegralLimzetaJay}.
\end{proof}
\begin{cor}
There is a measurable $\xi: [0,\infty)\rightarrow \Pee^\perp$ with finite $E^*$ variation such that $\xi(t)\in \bxi(t)$ and
$$
\|\bxi(t)\|_*=\|\xi(t)\|_*
$$
for all $t\ge 0$.
\end{cor}
\begin{proof}
We first claim that given any $\mu\in \Pee^\perp$, the function 
$$
f(t):=\inf\{\|\mu-\zeta\|_*: \zeta\in \bxi(t)\}\quad (t\ge 0)
$$
is Lebesgue measurable.  Note that we may rewrite this function as 
$$
f(t)=\inf\{\|\xi(t)-\mu+\eta\|_*: \eta\in \Cee^\perp\}=\|\bxi(t)-\bmu\|_*.
$$
By Corollary \ref{QuotientNormPeePerpLSC}, the norm on $\Pee^\perp/\Cee^\perp$ is lower-semicontinuous with convergence with respect to the metric 
induced by $\chi_\Cee^*$. Combined with the continuity of $\bxi$ \eqref{bxiCont}, $f$ is lower-semicontinuous for all but 
countably many times $t\ge 0$.  

\par By Lemma \ref{LSCfLem} in the appendix, $f$ is necessarily Lebesgue measurable. Since $\bxi(t)\subset \Pee^\perp$ is closed and nonempty for each $t\ge 0$,
it follows that 
$$
\bxi: [0,\infty)\rightsquigarrow \Pee 
$$
is a measurable, set-valued mapping (Theorem 8.1.4 of \cite{MR1048347}). As a result, the mapping 
$$
\bzeta: [0,\infty)\rightsquigarrow \Pee 
$$
defined by
$$
\bzeta(t):=\{\zeta\in \bxi(t): \|\zeta\|_*=\|\bxi(t)\|_*\}\quad (t\ge 0)
$$
is measurable (Theorem 8.2.11 of \cite{MR1048347}). Consequently, there is a measurable $\xi: [0,\infty)\rightarrow \Pee^\perp$ for which 
$$
\xi(t)\in \bzeta(t)
$$
and
$$
\|\xi(t)\|_*=\|\bxi(t)\|_*
$$
for all $t\ge 0$.  And since $E^*(\bxi,0,t)=E^*(\xi,0,t)$ for all $t\ge 0$, $\xi$ has finite $E^*$ variation. 
\end{proof}
\begin{rem}
We also note that as $\xi(0)\in \bxi(0)=\xi^0+\Cee^\perp$, $\xi$ satisfies the initial value condition \eqref{xizeroCond}.
\end{rem}
\begin{cor}
For almost every $t\ge 0$, 
$$
\lim_{j\rightarrow\infty}\|\zeta_{\tau_j}(t)\|_*=\|\xi(t)\|_*.
$$
\end{cor}
\begin{proof}
We have already noted that the limit 
$$
A(t)=\lim_{j\rightarrow\infty}\|\zeta_{\tau_j}(t)\|_*
$$
exists for all $t\ge 0$. In view of \eqref{IntegralLimzetaJay}, we may assume, after passing to a further subsequence if necessary, that 
$$
\lim_{j\rightarrow\infty}\chi^*_\Cee(\bzeta_{\tau_j}(t)-\bxi(t))=0
$$
for almost every $t\ge 0$. Recall that  
$$
\|\zeta_{\tau_j}(t)\|_*=\|\bzeta_{\tau_j}(t)\|_*
$$
by \eqref{chseekMin} for all but countable many $t\ge 0$. Therefore,
$$
A(t)=\lim_{j\rightarrow\infty}\|\bzeta_{\tau_j}(t)\|_*\ge\|\bxi(t)\|_*=\|\xi(t)\|_*
$$
for almost every $t\ge 0$. Here used Corollary \ref{QuotientNormPeePerpLSC}.

\par Since $\bg_{\tau}(t)\in \Jay^*(\xi_\tau(t))$ for almost all $t\ge 0$, we also have 
\begin{align}
\frac{1}{2}\|\xi(t)\|^2_*&\ge \frac{1}{2}\|\zeta_{\tau_j}(t)\|^2_*+\langle \xi(t)-\zeta_{\tau_j}(t),g_{\tau_j}(t)\rangle\\
&\ge  \frac{1}{2}\|\zeta_{\tau_j}(t)\|^2_*-\chi_\Cee^*(\bxi(t)-\bzeta_{\tau_j}(t)).
\end{align}
Thus, 
$$
\frac{1}{2}A(t)^2=\lim_{j\rightarrow\infty} \frac{1}{2}\|\zeta_{\tau_j}(t)\|^2_*\le  \frac{1}{2}\|\xi(t)\|^2_*
$$
for almost every $t\ge0$. We conclude that $A(t)=\|\xi(t)\|_*$ for almost every $t\ge 0$.
\end{proof}
\begin{rem}
An immediate corollary of the above assertion is that 
$$
\text{$[0,\infty)\ni t\mapsto\|\xi(t)\|_*$ is essentially nonincreasing.}
$$
\end{rem}

\subsection{Proof of the main theorem}
So far we have $\xi:[0,\infty)\rightarrow \Pee^\perp$, which is our candidate for the weak solution asserted to exist in Theorem \ref{mainThm}. We still need to establish 
$$
\bg(t)\in \Jay^*(\xi(t))
$$
for almost every $t\ge 0$, 
$$
[0,\infty)\ni t\mapsto E(\bg(t))\quad\text{is essentially nonincreasing},
$$
and 
$$
\frac{1}{2}\|\xi(s)\|_*^2- \frac{1}{2}\|\xi(t)\|_*^2= E^*V(\xi,s,t)+\int^t_sE(\bg(\tau))d\tau
$$
for almost every $0\le s\le t$. We'll start with verifying the following assertion. 
\begin{lem}\label{SubDiffLemmaGee}
For almost every $t\ge0$, $\bg(t)\in \Jay^*(\xi(t))$.
\end{lem}
\begin{proof}
For $\mu\in\Pee^\perp$ and $T>0$, 
\begin{align}
\int^T_0\frac{1}{2}\|\mu\|^2_*dt&\ge \int^T_0\left(\frac{1}{2}\|\xi_{\tau_j}(t)\|^2_*+\langle \mu-\xi_{\tau_j}(t),g_{\tau_j}(t)\rangle \right)dt \\
&=\int^T_0\frac{1}{2}\|\xi_{\tau_j}(t)\|^2_*dt +\int^T_0\langle \mu-\xi(t),g_{\tau_j}(t)\rangle dt +\int^T_0\langle \xi_{\tau_j}(t)-\xi(t),g_{\tau_j}(t)\rangle dt \\
&\ge \int^T_0\frac{1}{2}\|\xi_{\tau_j}(t)\|^2_*dt +\int^T_0\langle \mu-\xi(t),g_{\tau_j}(t)\rangle dt -\int^T_0 \chi_\Cee^*(\xi_{\tau_j}(t)-\xi(t)) dt \\
&= \int^T_0\frac{1}{2}\|\xi(t)\|^2_*dt +\int^T_0\langle \mu-\xi(t),g_{\tau_j}(t)\rangle dt +o(1)
\end{align}
as $j\rightarrow\infty$. Therefore, it suffices to verify 
\be\label{MainClaimSubdiffLemm}
\lim_{j\rightarrow\infty}\int^T_0\langle \mu-\xi(t),g_{\tau_j}(t)\rangle dt=\int^T_0\langle \mu-\xi(t),g(t)\rangle dt.
\ee

\par Note that $\rho(t):=\mu-\xi(t)\in \Pee^\perp$ satisfies
$$
\|\rho(t)\|_*\le \|\mu\|_*+\|\xi(t)\|_*\le \|\mu\|_*+\|\xi^0\|_*
$$ 
for almost every $t\ge0$. Proposition \ref{smoothlemma} in the appendix asserts that for each $\epsilon\in (0,1)$, there is 
$\rho^\epsilon: [0,\infty)\rightarrow (\Vee/\Pee)'$ such that: for all $t\ge 0$, $\rho^\epsilon(t)$ arises as the integration of a smooth function 
against $\sigma$ with
\be\label{rhoepsgteeEst}
|\langle \rho^\epsilon(t),g\rangle -\langle \rho(t),g\rangle|\le 4\sqrt{2}\epsilon\|\rho(t)\|_*
\ee
for each $g\in \Cee$. Moreover, $t\mapsto \langle \rho^\epsilon(t),h(t)\rangle$ is measurable for any measurable $\bh: [0,\infty)\rightarrow \Vee/\Pee$;
this is verified in Lemma \ref{measurablemappingRho}. Therefore,
$$
\left|\int^T_0\langle\rho^\epsilon(t),g_{\tau_j}(t)\rangle dt-\int^T_0\langle \rho(t),g_{\tau_j}(t)\rangle dt\right|
\le 4\sqrt{2\epsilon}T( \|\mu\|_*+\|\xi^0\|_*).
$$

\par We can also solve the PDE
\be\label{acheEqnRhoeps}
-(\Delta h+2h)=\rho^\epsilon(t)
\ee
weakly in $\Stwo$ for each $t\ge 0$ and $\epsilon \in (0,1)$.  For example, we can minimize the functional
$$
h\mapsto\int_{\Stwo}\left(\frac{1}{2}|\nabla h|^2-2h^2\right)d\sigma-\int_\Stwo \rho^\epsilon(t) hd\sigma
$$
uniquely among $h\in {\mathcal V}$ which satisfies 
$$
\int_\Stwo hu_id\sigma=0
$$
for $i=1,2,3$. The minimizing $h^\epsilon(t)\in {\mathcal V}$ satisfies 
$$
(\bh^\epsilon(t),{\bf f})=\int_\Stwo\left(\nabla h^\epsilon(t)\cdot \nabla f-2h^\epsilon(t)f\right)d\sigma=\int_{\Stwo}\rho^\epsilon(t) fd\sigma
$$
for all ${\bf f}\in {\mathcal V}/\Pee$.  Also note that Pettis' theorem (Chapter V section 4 of \cite{MR617913}) 
Lemma \ref{measurablemappingRho} imply $\bh^\epsilon:[0,\infty)\rightarrow {\mathcal V}/\Pee$ is measurable.

\par Observe 
\begin{align}
\int^T_0\langle\rho^\epsilon(t),g_{\tau_j}(t)\rangle dt&=\int^T_0\int_\Stwo\rho^\epsilon(t)g_{\tau_j}(t)d\sigma dt\\
&=\int^T_0\int_\Stwo\left(\nabla h^\epsilon(t)\cdot \nabla g_{\tau_j}(t)-2h^\epsilon(t)g_{\tau_j}(t)\right)d\sigma dt.
\end{align}
By the weak convergence of $(\bg_{\tau_j})_{j\in \N}$ to $\bg$ in $L^2([0,\infty);{\mathcal V}/\Pee)$,
\begin{align}
\lim_{j\rightarrow\infty}\int^T_0\langle\rho^\epsilon(t),g_{\tau_j}(t)\rangle dt&=
\int^T_0\int_\Stwo\left(\nabla h^\epsilon(t)\cdot \nabla g(t)-2h^\epsilon(t)g(t)\right)d\sigma dt\\
&=\int^T_0\int_\Stwo \rho^\epsilon(t)g(t)d\sigma dt\\
&=\int^T_0\langle\rho^\epsilon(t),g(t)\rangle dt.
\end{align}
Therefore,
\begin{align}
&\limsup_{j\rightarrow\infty}\left|\int^T_0\langle \rho(t),g_{\tau_j}(t)\rangle dt-\int^T_0\langle \rho(t),g(t)\rangle dt\right|\\
&\quad\le \limsup_{j\rightarrow\infty}\left|\int^T_0\langle \rho(t),g_{\tau_j}(t)\rangle dt-\int^T_0\langle \rho^\epsilon(t),g_{\tau_j}(t)\rangle dt\right|+\\
&\quad \quad  + \limsup_{j\rightarrow\infty}\left|\int^T_0\langle \rho^\epsilon(t),g_{\tau_j}(t)\rangle dt-\int^T_0\langle \rho(t),g(t)\rangle dt\right|\\
&\quad\le 4\sqrt{2\epsilon}T( \|\mu\|_*+\|\xi^0\|_*)+\left|\int^T_0\langle \rho^\epsilon(t),g(t)\rangle dt-\int^T_0\langle \rho(t),g(t)\rangle dt\right|.
\end{align}
Again recalling that $g(t)\in \Cee$ for all $t\ge 0$ and using inequality  \eqref{rhoepsgteeEst} leads to
$$
\limsup_{j\rightarrow\infty}\left|\int^T_0\langle \rho(t),g_{\tau_j}(t)\rangle dt-\int^T_0\langle \rho(t),g(t)\rangle dt\right|\le 8\sqrt{2\epsilon}T( \|\mu\|_*+\|\xi^0\|_*).
$$
Since $\epsilon\in (0,1)$ was arbitrary, we conclude \eqref{MainClaimSubdiffLemm}.
\end{proof}
Next we claim that an energy inequality holds. 
\begin{lem}
For almost every $0\le s\le t$, 
$$
\frac{1}{2}\|\xi(s)\|_*^2- \frac{1}{2}\|\xi(t)\|_*^2\ge E^*V(\xi,s,t)+\int^t_sE(\bg(\tau))d\tau.
$$
\end{lem}
\begin{proof}
Assume $0<s<t$ and choose $0<\delta<s$. Further select $\tau<\min\{\delta,t-s\}$.  Since $s+\tau<t$, we may select $j,k\in \N$ with $j<k$ such that $(j-1)\tau<s\le j\tau$ and $(k-1)\tau<t\le k\tau$. Likewise, since $\tau<\delta$, $s-\delta<(j-1)\tau$; so there is $i\le  j-1$ for which $(i-1)\tau<s-\delta\le i\tau $.
By \eqref{DiscreteMonForm1},
\begin{align}
\frac{1}{2}\|\zeta_{\tau}(s-\delta)\|_*^2&=\frac{1}{2}\|\xi^i\|_*^2\\
&\ge \frac{1}{2}\|\xi^{j-1}\|_*^2\\
&\ge \frac{1}{2}\|\xi^{k}\|_*^2+\sum^{k}_{\ell=j}\tau \left[E^*\left(\frac{\xi^{\ell}-\xi^{\ell-1}}{\tau}\right)+E(\bg^\ell)\right]\\
&=\frac{1}{2}\|\zeta_{\tau}(t)\|_*^2+E^*V(\xi_\tau,(j-1)\tau,k\tau)+\int^{k\tau}_{(j-1)\tau}E(\bg_\tau(r))dr\\
&\ge\frac{1}{2}\|\zeta_{\tau}(t)\|_*^2+E^*V(\xi_\tau,s,t)+\int^{t}_{s}E(\bg_\tau(r))dr.
\end{align}
We can then pass to the limit along an appropriate sequence of $\tau$ tending to $0$ to find 
$$
A(s-\delta)\ge A(t)+E^*V(\xi,s,t)+\int^{t}_{s}E(\bg(r))dr.
$$ 
\par Since $A$ is monotone, it is continuous except for possibly on a countable set of times. Therefore, we 
may send $\delta$ to $0$ and find 
$$
A(s)\ge A(t)+E^*V(\xi,s,t)+\int^{t}_{s}E(\bg(r))dr
$$
for all but countably many times $0\le s\le t$. We conclude upon recalling that $A(t)=\frac{1}{2}\|\xi(t)\|^2_*$ for almost every time $t\ge 0$.
\end{proof}

\par Now we are in position to establish the monotonicity of $E\circ \bg$. 
\begin{lem}
The function $E\circ \bg$ is essentially nonincreasing.
\end{lem}
\begin{proof}
By the previous lemma, 
\be\label{EcircbgNonINeq1}
\frac{1}{2}\|\xi(s)\|_*^2- \frac{1}{2}\|\xi(t)\|_*^2\ge E^*V(\xi,s,t)+\int^t_sE(\bg(\tau))d\tau
\ee
for almost all $0\le s\le t$. We also have that $\bg(s)\in \Jay^*(\xi(s))$ and $g(s)\in \Cee$ for almost every $s\ge 0$. For $t$ larger than such an $s$,
\begin{align}\label{SubdiffINequalityTandS}
\frac{1}{2}\|\xi(t)\|_*^2&\ge \frac{1}{2}\|\xi(s)\|_*^2+\langle \xi(t)-\xi(s),g(s)\rangle\nonumber\\
&\ge \frac{1}{2}\|\xi(s)\|_*^2+(t-s)\left\langle \frac{\xi(t)-\xi(s)}{t-s},g(s)\right\rangle\nonumber\\
&\ge \frac{1}{2}\|\xi(s)\|_*^2-(t-s)E^*\left(\frac{\xi(t)-\xi(s)}{t-s}\right)-(t-s)E(\bg(s))\nonumber\\
&\ge \frac{1}{2}\|\xi(s)\|_*^2-E^*V(\xi,s,t)-(t-s)E(\bg(s)).
\end{align}
That is, 
\be\label{EcircbgNonINeq2}
\frac{1}{2}\|\xi(s)\|_*^2- \frac{1}{2}\|\xi(t)\|_*^2\le E^*V(\xi,s,t)+(t-s)E(\bg(s)).
\ee

\par Comparing \eqref{EcircbgNonINeq1} and \eqref{EcircbgNonINeq2}, we find 
$$
\int^t_sE(\bg(\tau))d\tau\le(t-s)E(\bg(s))
$$
for almost all $0\le s\le t$. By Lemma \ref{NonincreasingCondAppendix} in the appendix, 
$E\circ \bg$ is necessarily essentially nonincreasing. 
\end{proof}
The final detail needed in our proof of Theorem \ref{mainThm} is that equality (essentially) holds in the energy identity. 

\begin{lem}
For almost every $s\le t$, 
$$
\frac{1}{2}\|\xi(s)\|_*^2- \frac{1}{2}\|\xi(t)\|_*^2= E^*V(\xi,s,t)+\int^t_sE(\bg(\tau))d\tau.
$$
\end{lem}
\begin{proof}
By Lemma \ref{EssNonIncLemma}, we may choose a  nondecreasing function $I:[0,\infty)\rightarrow [0,\infty)$ such that 
$$
I(t)=E(\bg(t))
$$
for almost every $t\ge 0$. Let $0\le s< t$ be such that $I(s)=E(\bg(s))$, $I(t)=E(\bg(t))$, 
and $\bg(s)\in \Jay^*(\xi(s))$. Since $I$ is nonincreasing and $0\le I(\tau)\le E(\bg^0)$ for almost every $\tau\ge 0$, $I$ is Riemann integrable on $[s,t]$.  

\par  Fix $\epsilon>0$. We may select $\delta>0$ such that for any partition $s=t_0<t_1<\dots<t_N=t$ with $\max_{1\le i\le N}(t_i-t_{i-1})\le\delta$ implies
$$
\sum^N_{i=1}(t_i-t_{i-1})(I(t_{i-1})-I(t_i))\le \epsilon. 
$$
In view of Lemma \ref{SubDiffLemmaGee}, we may choose such a partition for which $\bg(t_i)\in \Jay^*(\xi(t_i))$ and $I(t_i)=E(\bg(t_i))$ for $i=1,\dots, N-1.$ Employing inequality \eqref{SubdiffINequalityTandS}, 
we find 
\begin{align}
\frac{1}{2}\|\xi(s)\|_*^2- \frac{1}{2}\|\xi(t)\|_*^2&=\sum^N_{i=1}\left[\frac{1}{2}\|\xi(t_{i-1})\|_*^2- \frac{1}{2}\|\xi(t_i)\|_*^2\right]\\
&\le \sum^N_{i=1}\left[ E^*V(\xi,t_{i-1},t_i)+(t_i-t_{i-1})I(t_{i-1})\right]\\
&= E^*V(\xi,s,t)+ \sum^N_{i=1}(t_i-t_{i-1})I(t_{i-1})\\
&\le E^*V(\xi,s,t)+ \sum^N_{i=1}(t_i-t_{i-1})I(t_{i})+\epsilon\\
&\le E^*V(\xi,s,t)+ \int^t_sI(\tau)d\tau+\epsilon\\
&= E^*V(\xi,s,t)+ \int^t_sE(\bg(\tau))d\tau+\epsilon.
\end{align} 
\par Since $\epsilon$ was arbitrary, 
$$
\frac{1}{2}\|\xi(s)\|_*^2- \frac{1}{2}\|\xi(t)\|_*^2\le E^*V(\xi,s,t)+ \int^t_sE(\bg(\tau))d\tau.
$$
We conclude by noting that this inequality and the opposite inequality holds for almost every $0\le s\le t$.
\end{proof}

\thanks {\bf Acknowledgements}:  This research was partially supported by NSF award DMS-1554130.  Figures \ref{FirstMeissnerFigure} and \ref{SecondMeissnerFigure} were made with the Meissner Tetrahedra -- Wolfram Demonstration Project.  We would also like to thank Kristine Kreidler and Cheikh N'diaye for their interest in this work and for their encouragement.

\appendix

\section{Smoothing elements of $\Pee^\perp$}
We will discuss a way to smooth measures on $\Stwo$ which is tailored for the needs of this paper. With this goal in mind, we will choose a family $(\psi^\epsilon)_{\epsilon\in (0,1)}$ satisfying 
$$
\begin{cases}
\psi^\epsilon\in C^\infty(\R)\\
\psi^\epsilon\ge 0\\
\text{supp}(\psi^\epsilon)=[1-\epsilon,1]\\
\end{cases}
$$
and 
$$
\int_{\Stwo}\psi^\epsilon(u_1)d\sigma(u)=1
$$
for each $\epsilon\in (0,1)$. Since $\sigma$ is invariant under orthogonal transformations, 
$$
\int_{\Stwo}\psi^\epsilon(u\cdot v)d\sigma(u)=1
$$
for each $v\in \Stwo$. Moreover, the support of $u\mapsto \psi^\epsilon(u\cdot v)$ is the spherical cap $\{u\in \Stwo: u\cdot v\ge 1-\epsilon\}$.

\par For a given $\mu\in M(\Stwo)$, we can define 
\be\label{tildeMuEps}
\tilde\mu^\epsilon(u):=\int_{\Stwo}\psi^\epsilon(u\cdot v)d\mu(v)
\ee
for $u\in \Stwo$. It is routine to check that $\tilde\mu^\epsilon\in C^\infty(\Stwo)$.  We will identify this function with the measure it induces when 
integrated against $\sigma$:
$$
\int_\Stwo gd\tilde\mu^\epsilon:=\int_{\Stwo}g(u)\tilde\mu^\epsilon(u)d\sigma(u)=\int_{\Stwo}\int_{\Stwo}g(u)\psi^\epsilon(u\cdot v)d\mu(v)d\sigma(u)
$$
for $g\in C(\Stwo)$. 

\begin{lem}
Suppose $\mu\in M(\Stwo)$ and $g\in C(\Stwo)$ satisfies
$$
|g(u)-g(v)\le |u-v|
$$
for $u,v\in \Stwo$. Then 
\be\label{SmoothingEstLipGee}
\left|\int_\Stwo gd\tilde\mu^\epsilon-\int_\Stwo gd\mu\right|\le  \sqrt{2\epsilon}\;\|\mu\|_*.
\ee
\end{lem}
\begin{proof}
Observe
\begin{align}
\int_\Stwo gd\tilde\mu^\epsilon-\int_\Stwo gd\mu&=\int_{\Stwo}\int_{\Stwo}g(u)\psi^\epsilon(u\cdot v)d\mu(v)d\sigma(u)-\int_\Stwo g(v)d\mu(v)\\
&=\int_{\Stwo}\left[\int_{\Stwo}g(u)\psi^\epsilon(u\cdot v)d\sigma(u)- g(v)\right]d\mu(v)\\
&=\int_{\Stwo}\left[\int_{\Stwo}(g(u)-g(v))\psi^\epsilon(u\cdot v)d\sigma(u)\right]d\mu(v)\\
&=\int_{\Stwo}\left[\int_{\{u\cdot v\ge 1-\epsilon\}}(g(u)-g(v))\psi^\epsilon(u\cdot v)d\sigma(u)\right]d\mu(v)\\
&\le \int_{\Stwo}\left[\int_{\{u\cdot v\ge 1-\epsilon\}}|g(u)-g(v)|\psi^\epsilon(u\cdot v)d\sigma(u)\right]d|\mu|(v)\\
&\le \int_{\Stwo}\left[\int_{\{u\cdot v\ge 1-\epsilon\}}|u-v|\psi^\epsilon(u\cdot v)d\sigma(u)\right]d|\mu|(v).
\end{align}
Note that on the spherical cap $\{u\in \Stwo: u\cdot v\ge 1-\epsilon\}$,
$$
|u-v|^2=2(1-u\cdot v)\le 2\epsilon.
$$
As a result, 
$$
\int_\Stwo gd\tilde\mu^\epsilon-\int_\Stwo gd\mu\le \sqrt{2\epsilon}\int_{\Stwo}\left[\int_{\Stwo}\psi^\epsilon(u\cdot v)d\sigma(u)\right]d|\mu|(v)=\| \mu\|_*.
$$
We can use a similar argument to obtain the same upper bound for $\int_\Stwo gd\mu-\int_\Stwo gd\tilde\mu^\epsilon$.
\end{proof}
We will now introduce the approximation
\be\label{myooEps}
\mu^\epsilon(u):=\tilde\mu^\epsilon(u)-a^\epsilon\cdot u,
\ee
where $a^\epsilon\in \R^3$ is chosen so that $\mu^\epsilon$ is orthogonal to $\Pee$ in $L^2(\Stwo)$. That is, 
$$
a^\epsilon_i=\frac{\displaystyle\int_{\Stwo}u_i\tilde\mu^\epsilon(u)d\sigma(u)}{\displaystyle\int_{\Stwo}u_i^2d\sigma(u)}, \quad i=1,2,3.
$$
\begin{prop}\label{smoothlemma}
Suppose $g\in \Cee$, $\mu\in \Pee^\perp$, and define $\mu^\epsilon$ by \eqref{myooEps}. Then 
$\mu^\epsilon\in C^\infty(\Stwo)$, as a measure $\mu^\epsilon\in\Pee^\perp$, and
\be
\left|\int_\Stwo gd\mu^\epsilon-\int_\Stwo gd\mu\right|\le  4\sqrt{2\epsilon}\;\|\mu\|_*.
\ee
\end{prop}
\begin{proof}
It is clear that $\mu^\epsilon\in C^\infty(\Stwo)$. Since $\mu^\epsilon$ is orthogonal to $\Pee$ in $L^2(\Stwo)$, the measure it induces belongs to 
$\Pee^\perp$  as
$$
\int_{\Stwo}u_id\mu^\epsilon(u)=\int_{\Stwo}u_i\mu^\epsilon(u)d\sigma(u)=0
$$
($i=1,2,3$).

\par Now let $H$ be the support function of a constant width body $K$ associated with $g$. Let $a\in K$ be the center of the circumball of $K$.
 By Corollary \ref{EstimatesSuppHnoted},
$$
\begin{cases}
|H(u)-1/2-a\cdot u|\le 1
\\
|H(u)-H(v)-a\cdot (u-v)|\le |u-v|
\end{cases}
$$
for all $u,v\in \Stwo$.  Therefore, if we set $\tilde g(u)=g(u)-a\cdot u$, 
$$
\begin{cases}
|\tilde g(u)|\le 1\\
|\tilde g(u)-\tilde g(v)|\le |u-v|
\end{cases}
$$
for all $u,v\in \Stwo$. 
\par By inequality \eqref{SmoothingEstLipGee},
\begin{align}
\left|\int_\Stwo gd\mu^\epsilon-\int_\Stwo gd\mu\right|
&\quad = \left|\int_\Stwo \tilde gd\mu^\epsilon-\int_\Stwo \tilde gd\mu\right|\\
&\quad = \left|\int_\Stwo \tilde gd\tilde\mu^\epsilon-\int_\Stwo \tilde gd\mu-\int_\Stwo \tilde g(u) a^\epsilon\cdot u d\sigma(u)\right|\\
&\quad \le \sqrt{2\epsilon}\|\mu\|_* +\left|\int_\Stwo \tilde g(u) a^\epsilon\cdot u d\sigma(u)\right|\\
&\quad \le \sqrt{2\epsilon}\|\mu\|_* +4\pi |a^\epsilon|.
\end{align}
Recall  
$$
\int_\Stwo u_i^2d\sigma(u)=\frac{4\pi}{3}
$$
for $i=1,2,3$ and $\mu\in \Pee^\perp.$ It follows that  
\begin{align}
\frac{4\pi}{3} a^\epsilon&=\int_{\Stwo}u\tilde\mu^\epsilon(u)d\sigma(u)\\
&=\int_{\Stwo}\int_{\Stwo}uk^\epsilon(u\cdot v)d\mu(v)d\sigma(u)\\
&=\int_{\Stwo}\int_{\Stwo}(u-v)k^\epsilon(u\cdot v)d\sigma(u)d\mu(v) +\int_{\Stwo}\int_{\Stwo}vk^\epsilon(u\cdot v)d\mu(v)d\sigma(u)\\
&=\int_{\Stwo}\int_{\Stwo}(u-v)k^\epsilon(u\cdot v)d\sigma(u)d\mu(v) +\int_{\Stwo}v\left(\int_{\Stwo}k^\epsilon(u\cdot v)d\sigma(u)\right)d\mu(v)\\
&=\int_{\Stwo}\int_{\Stwo}(u-v)k^\epsilon(u\cdot v)d\sigma(u)d\mu(v)+\int_\Stwo vd\mu(v)\\
&=\int_{\Stwo}\int_{\Stwo}(u-v)k^\epsilon(u\cdot v)d\sigma(u)d\mu(v).
\end{align}
Therefore, 
\begin{align}
4\pi |a^\epsilon|&\le 3\int_{\Stwo}\int_{\Stwo}|u-v|k^\epsilon(u\cdot v)d\sigma(u)d|\mu|(v)\\
&= 3\int_{\Stwo}\int_{u\cdot v\ge 1-\epsilon}|u-v|k^\epsilon(u\cdot v)d\sigma(u)d|\mu|(v)\\
&\le 3\sqrt{2\epsilon}\int_{\Stwo}\int_{\Stwo}k^\epsilon(u\cdot v)d\sigma(u)d|\mu|(v)\\
&=3\sqrt{2\epsilon}\|\mu\|_*.
\end{align}
Consequently, 
$$
\left|\int_\Stwo gd\mu^\epsilon-\int_\Stwo gd\mu\right|\le \sqrt{2\epsilon}\|\mu\|_*+3\sqrt{2\epsilon}\|\mu\|_*.
$$
\end{proof}

\par We will need one more technical assertion regarding how this smoothing technique is applied to measurable mappings $\rho:[0,\infty)\rightarrow \Pee^\perp$. We will employ the fact that if $\mu\in\Pee^\perp$, then 
\be\label{VeeOverPeeBound}
\left|\int_{\Stwo}\mu^\epsilon h d\sigma\right| \le \frac{1}{2}\left(\int_{\Stwo}(\mu^\epsilon)^2d\sigma\right)^{1/2}(\bh,\bh)^{1/2}
\ee
for all $\bh\in \Vee/\Pee$; we established this inequality in our proof of Lemma \ref{AchOneMeaswithCeegivesMeas}. In particular, 
$$
{\mathcal V}/\Pee \ni \bh\mapsto \int_{\Stwo} h\mu^\epsilon d\sigma 
$$
is continuous, so we may consider $\mu^\epsilon$ an element of $(\Vee/\Pee)'$.
\begin{lem}\label{measurablemappingRho}
Suppose $\rho:[0,\infty)\rightarrow \Pee^\perp$ is measurable and define $\rho^\epsilon:[0,\infty)\rightarrow (\Vee/\Pee)'$ by 
$$
\rho^\epsilon(t):=(\rho(t))^\epsilon
$$
for each $t\ge 0$.  Then for all  $\bg \in {\mathcal V}/\Pee$,
$$
[0,\infty)\ni t\mapsto\int_{\Stwo} h\rho^\epsilon(t) d\sigma 
$$
is measurable. Moreover, for all measurable $\bg: [0,\infty)\rightarrow {\mathcal V}/\Pee$,
$$
[0,\infty)\ni t\mapsto\int_{\Stwo} g(t)\rho^\epsilon(t) d\sigma 
$$
is measurable.
\end{lem}
\begin{proof}
Let $\bh \in \Vee/\Pee$.  For a given $t\ge 0$, set $\mu=\rho(t)$ and define $\tilde\rho^\epsilon(t):=\tilde\mu^\epsilon$ as in \eqref{tildeMuEps}.  Note that 
$$
\int_{\Stwo}h \tilde\rho^\epsilon(t)d\sigma =\langle \rho(t),h^\epsilon\rangle,
$$
where 
$$
h^\epsilon(u):=\int_{\Stwo}\psi^\epsilon(u\cdot v)g(v)d\sigma(v)
$$
for $u\in\Stwo$.  We also have 
\begin{align}
\int_{\Stwo}\rho^\epsilon(t)hd\sigma & = \int_{\Stwo}h\tilde\rho^\epsilon(t)d\sigma -\int_\Stwo a^\epsilon(t)\cdot uh(u)d\sigma(u)\\
& =\langle \rho(t),h^\epsilon\rangle -\sum^3_{i=1}a^\epsilon_i(t) \left(\int_\Stwo u_ih(u)d\sigma(u)\right),
\end{align}
where 
$$
a_i^\epsilon(t)=\frac{\displaystyle\int_{\Stwo}u_i\tilde\rho^\epsilon(t)d\sigma}{\displaystyle\int_{\Stwo}(u_i)^2d\sigma}=\langle \rho^\epsilon(t),\phi_i\rangle 
=\langle \rho(t),\phi^\epsilon_i\rangle
$$
for $i=1,2,3$. Here $\phi_i(v)=v_i/\int_{\Stwo}(u_i)^2d\sigma$.  Thus, 
$$
\int_{\Stwo}\rho^\epsilon(t)gd\sigma =\langle \rho(t),g^\epsilon\rangle-\sum^3_{i=1} \left(\int_\Stwo u_ig(u)d\sigma(u)\right)\langle \rho(t),\phi^\epsilon_i\rangle
$$
is a measurable function since $\rho$ is measurable.  

\par Observe that for any fixed $t$,
$$
\Vee/\Pee\ni \bh \mapsto\int_{\Stwo} h\rho^\epsilon(t) d\sigma 
$$
is continuous. This continuity follows directly from inequality \eqref{VeeOverPeeBound}. Therefore, 
$$
[0,\infty)\times\Vee/\Pee\ni  (t,\bh) \mapsto\int_{\Stwo} h\rho^\epsilon(t) d\sigma
$$
is a Carath\'eodory function. As a result, if $\bg: [0,\infty)\rightarrow {\mathcal V}/\Pee$ is measurable, so is 
$$
[0,\infty)\ni t\mapsto\int_{\Stwo} g(t)\rho^\epsilon(t) d\sigma 
$$
(Lemma 8.2.3 of \cite{MR1048347}).
\end{proof}

\section{Essentially monotone functions}\label{EssentialMonotoneSect}
We'll say that a function $f:[0,\infty)\rightarrow \R$ is {\it essentially nondecreasing} if there is a set $N\subset [0,\infty)$ of measure 0 for which $f(t)\ge f(s)$ whenever $t\ge s$ and $t,s\not\in N$. If under the same conditions, $f(t)\le f(s)$,  we'll say $f$ is {\it essentially nonincreasing}. If $f$ is either essentially nondecreasing or nonincreasing, we'll say $f$ is {\it essentially monotone}. 

\par  It turns out that the only way essentially monotone functions arise is by altering monotone functions on sets of measure $0$. 
\begin{lem}\label{EssNonIncLemma}
Suppose $f:[0,\infty)\rightarrow \R$ is essentially nondecreasing. There is a nondecreasing function $g: [0,\infty)\rightarrow [-\infty,\infty)$ such that $f(t)=g(t)$ for almost every $t\ge 0$. 
\end{lem}

\begin{proof}
Choose a null set $N\subset [0,\infty)$ such that $f|_{N^c}$ is nondecreasing. For $t\ge 0$, let us also define
$$
g(t):=\sup\left\{c\in \R: c\le f(\tau),\;\tau\ge t,\tau\not\in N\right\}.
$$  
Suppose $t \not\in N$. Then $c\le f(t)$ in the definition above; therefore, $g(t)\le f(t)$. And as $f(t)\le f(\tau)$ for $\tau\ge t$ and $\tau\not\in N$, we also have $ f(t)\le g(t)$.  Thus, $f$ and $g$ agree on $N^c$.  

\par Let us check that $g$ is nondecreasing. Suppose $t_1\le t_2$. Observe that if $c\le f(\tau)$ for all $\tau\ge t_1$ with $\tau\not\in N$, then $c\le f(\tau)$ for all $\tau\ge t_2$ with $\tau\not\in N$. In this case, $c\le g(t_2)$. Therefore, $g(t_1)\le g(t_2)$. Also note that for any $t\ge 0$, there is $\tau>t$ with $\tau\not\in N$. Therefore, $g(t)\le g(\tau)=f(\tau)<\infty$. 
\end{proof}
We now generalize an elementary fact: a Riemann integrable function from $[0,\infty)$ into $[0,\infty)$ which is nonincreasing tends to $0$ faster than $1/t$ as $t\rightarrow\infty$. 
\begin{prop}\label{teeffteeGoesZero}
Suppose that $f: [0,\infty)\rightarrow [0,\infty)$ is essentially nonincreasing and 
$$
\int^\infty_0f(t)dt<\infty.
$$
Then there is a null set $N\subset [0,\infty)$ for which 
\be\label{oneovertexercise}
\lim_{\substack{t\rightarrow\infty \\ t\not\in N}}tf(t)=0.
\ee
\end{prop}
\begin{proof}
Let $N\subset [0,\infty)$ be a null set for which $f|_{N^c}$ is nonincreasing. Next, choose an unbounded, increasing sequence of times $(t_k)_{k\in \N}\subset N^c$ and fix $\epsilon>0$. Since $f$ is integrable, there is $K(\epsilon)\in\N$ such that 
$$
\frac{t_k}{2}f(t_k)\le \int^{t_k}_{t_k/2}f(t)dt\le \frac{1}{2}\epsilon
$$
for $k\ge K(\epsilon)$. That is, $t_kf(t_k)\le\epsilon$ for $k\ge K(\epsilon)$. Therefore, $\lim_{k\rightarrow\infty}t_kf(t_k)=0$. We conclude \eqref{oneovertexercise}.
\end{proof}

\par Note that if $f:[0,\infty)\rightarrow \R$ is essentially nonincreasing, then 
\be\label{WeakNonincreasingCondition}
\int^t_sf(\tau)d\tau\le (t-s)f(s)
\ee
for almost every $0\le s\le t<\infty$. The converse is also true. 

\begin{lem}\label{NonincreasingCondAppendix}
Suppose $f:[0,\infty)\rightarrow \R$ is locally integrable and 
 \eqref{WeakNonincreasingCondition} holds for almost every $0\le s\le t<\infty$. Then $f$ is essentially nonincreasing.
\end{lem}
\begin{proof}
1. First suppose that $f$ is smooth in $(0,\infty)$.  This implies $f$ satisfies \eqref{WeakNonincreasingCondition} for all $0<s\le t<\infty$. We can integrate by parts to find 
\begin{align}
\int^t_sf(\tau)d\tau &= \int^t_sf(\tau)\frac{d}{d\tau}(\tau-t)d\tau\\
&= (\tau-t)f(\tau)\Big|^t_s-\int^t_s f'(\tau)(\tau-t)d\tau \\
&=f(s)(t-s)+\int^t_s f'(\tau)(t-\tau)d\tau\\
&\le f(s)(t-s).
\end{align}
Consequently, 
\be\label{WeakNonincreasingCondition2}
\int^t_s f'(\tau)(t-\tau)d\tau\le 0
\ee
for all $0<s\le t<\infty$. 

\par We claim that $f'(\tau)\le 0$ for all $0<\tau<\infty$. Otherwise, $f'(s_0)> 0$ for some $s_0\in (0,\infty)$. By continuity, there is $\delta>0$ for which 
$$
f'(\tau)\ge \frac{1}{2}f'(s_0)
$$
for $s_0\le \tau<s_0+\delta$. In particular, 
$$
\int^{s_0+\delta}_{s_0} f'(\tau)((s_0+\delta)-\tau)d\tau\ge
 \frac{1}{2}f'(s_0)\cdot\int^{s_0+\delta}_{s_0}((s_0+\delta)-\tau)d\tau= \frac{1}{2}f'(s_0)\cdot \frac{\delta^2}{2}>0.
$$
This would contradict \eqref{WeakNonincreasingCondition2} with $s=s_0$ and $t=s_0+\delta$. As a result, if $f$ is smooth and satisfies \eqref{WeakNonincreasingCondition}, then $f$ is necessarily nonincreasing. 

\par 2. Now let us suppose $f$ and $f^\epsilon=\eta^\epsilon*f$ is the standard mollification of $f$ which is defined on $(\epsilon,\infty)$ (Chapter 4 of \cite{MR3409135}).  Let $\epsilon<s<t<\infty$. Observe that 
\begin{align}
\int^t_sf^\epsilon(\tau)d\tau&=\int^t_s\left(\int^\epsilon_{-\epsilon}\eta^\epsilon(r)f(\tau-r)dr\right)d\tau\\
&=\int^\epsilon_{-\epsilon}\eta^\epsilon(r)\left(\int^t_sf(\tau-r)d\tau\right)dr\\
&=\int^\epsilon_{-\epsilon}\eta^\epsilon(r)\left(\int^{t-r}_{s-r}f(\rho)d\rho\right)dr.
\end{align}
As \eqref{WeakNonincreasingCondition} holds outside of a null set, 
$$
\int^{t-r}_{s-r}f(\rho)d\rho\le ((t-r)-(s-t))f(s-r)=(t-s)f(s-r)
$$
for almost every $|r|<\epsilon$. Therefore, 
$$
\int^t_sf^\epsilon(\tau)d\tau\le (t-s)\int^\epsilon_{-\epsilon}\eta^\epsilon(r)f(s-r)dr=(t-s)f^\epsilon(s).
$$
By part 1 of this proof, $f^\epsilon$ is nonincreasing on $(\epsilon,\infty)$.  Since $f^\epsilon(s)\rightarrow f(s)$ for almost every $s\in (0,\infty)$, we have that $f$ is essentially nonincreasing. 
\end{proof}

\section{Essentially lower-semicontinuous functions}
We will say that a function $f:[0,\infty)\rightarrow \R$ is {\it essentially lower-semicontinuous} if $f$ is lower-semicontinuous at $t$ for almost every $t\ge 0$. A typical lower-semicontinous 
function $g: [0,\infty)\rightarrow \R$ is Borel measurable as its sublevel set $g^{-1}((-\infty,y])$ is closed in $[0,\infty)$ for each $y\in \R$.  It turns out that an essentially lower-semicontinuous function is at least Lebesgue measurable. 

\begin{prop}\label{LSCfLem}
Suppose $f:[0,\infty)\rightarrow \R$ is essentially lower-semicontinuous. Then $f$ is Lebesgue measurable. 
\end{prop}
\begin{proof}
Let $N\subset[0,\infty)$ be a null set for which $f$ may not be lower-semicontinuous at times $t\in N$. For a given $y\in \R$, we claim that 
\be\label{finvMeasurabilityClaim}
f^{-1}((y,\infty))\cap N^c=O\cap N^c
\ee
for some open $O\subset \R$.   If we can establish this claim, then 
\begin{align}
f^{-1}((y,\infty))&=f^{-1}((y,\infty))\cap (N\cup N^c)\\
&=\left(f^{-1}((y,\infty))\cap N\right)\cup \left(f^{-1}((y,\infty))\cap N^c\right)\\
&=\left(f^{-1}((y,\infty))\cap N\right)\cup \left(O\cap N^c\right)\\
&=\left(f^{-1}((y,\infty))\cap N\right)\cup \left(O^c\cup N\right)^c.
\end{align}
Since $f^{-1}((y,\infty))\cap N$ is a null set and $O^c\cup N$ is Lebesgue measurable, $f^{-1}((y,\infty))$ would be Lebesgue measurable. 

\par Let us now verify \eqref{finvMeasurabilityClaim}. Let $t\in f^{-1}((y,\infty))\cap N^c$. That is, $f(t)>y$ and $f$ is lower-semicontinuous at $t$. We note $f(t)-\epsilon>y$ for some 
$\epsilon>0$. We may also select $\delta_t>0$ such that if $|t-s|< \delta_t$ and $s\ge 0$, then
$$
f(t)\le f(s)+\epsilon.
$$
It follows that $f(s)\ge f(t)-\epsilon>y$ and thus
\be\label{relativelyOpenInclusion}
(t-\delta_t,t+\delta_t)\cap N^c\subset f^{-1}((y,\infty))\cap N^c.
\ee
Let us choose 
$$
O:=\bigcup\{(t-\delta_t,t+\delta_t): t\in f^{-1}((y,\infty))\cap N^c \}
$$
and note \eqref{relativelyOpenInclusion} implies 
$$
O\cap N^c\subset f^{-1}((y,\infty))\cap N^c.
$$
Alternatively, if $t\in f^{-1}((y,\infty))\cap N^c$, then $t\in (t-\delta_t,t+\delta_t)\cap N^c\subset O\cap N^c$. We conclude  $f^{-1}((y,\infty))\cap N^c\subset O\cap N^c$. This verifies \eqref{finvMeasurabilityClaim}.
\end{proof}

\appendix

\bibliography{ConstantWidthFlowArticleBib}{}

\bibliographystyle{amsplain}

\typeout{get arXiv to do 4 passes: Label(s) may have changed. Rerun}

\end{document}